\PassOptionsToPackage{table}{xcolor}
\documentclass[review,onefignum,onetabnum,nohypdvips]{siamart171218}

\usepackage{mathrsfs,amsmath,amssymb,bm}
\usepackage{lipsum}
\usepackage{amsfonts}
\usepackage{caption, graphicx, subfigure}
\usepackage{epstopdf}
\usepackage{makecell, rotating, bbding}
\usepackage{multirow}
\usepackage{algorithmic, algorithm}
\usepackage{enumerate}
\usepackage{booktabs}
\usepackage[misc]{ifsym}
\usepackage{mathtools}
\usepackage{amsmath}
\newtheorem{thm}{Theorem}[section]

\newtheorem{remark}[thm]{Remark}
\newtheorem{example}[thm]{Example}

\usepackage[draft]{changes}
\usepackage{lipsum}% <- For dummy text

\newsiamremark{hypothesis}{Hypothesis}
\crefname{hypothesis}{Hypothesis}{Hypotheses}
\newsiamthm{claim}{Claim}

\usepackage{amsopn}

\makeatletter
\newcommand*{\addFileDependency}[1]{% argument=file name and extension
  \typeout{(#1)}% latexmk will find this if $recorder=0 (however, in that case, it will ignore #1 if it is a .aux or .pdf file etc and it exists! if it doesn't exist, it will appear in the list of dependents regardless)
  \@addtofilelist{#1}% if you want it to appear in \listfiles, not really necessary and latexmk doesn't use this
  \IfFileExists{#1}{}{\typeout{No file #1.}}% latexmk will find this message if #1 doesn't exist (yet)
}
\makeatother

\headers{Accurately recover global quasiperiodic systems by finite points}{K. Jiang, Q. Zhou, P. Zhang}

\title{Accurately recover global quasiperiodic systems by finite points
\thanks{Submitted to the editors DATE.
}}

\author{Kai Jiang\thanks{
		Hunan Key Laboratory for Computation and Simulation in Science and Engineering,
		Key Laboratory of Intelligent Computing and Information Processing of Ministry
		of Education, School of
		Mathematics and Computational Science, Xiangtan University, Xiangtan, Hunan,
		China, 411105.
		(\email{kaijiang@xtu.edu.cn}, \email{qizhou@smail.xtu.edu.cn}).}
	\and Qi Zhou\footnotemark[2]
	\and Pingwen Zhang\thanks{
		School of Mathematics and Statistics, Wuhan University, Wuhan 430072, China. 
		School of Mathematical Sciences,
		Peking University, Beijing 100871, China.
		(\email{pzhang@pku.edu.cn}). } 
	}

\makeatletter

\newcommand{\Rmnum}[1]{\expandafter\@slowromancap\romannumeral #1@}
\makeatother

\begin{document}

\maketitle

\begin{abstract}
Quasiperiodic systems, related to irrational numbers, are space-filling structures without decay nor translation invariance. How to accurately recover these systems, especially for non-smooth cases, presents a big challenge in numerical computation.  
In this paper, we propose a new algorithm, finite points recovery (FPR) method, which is available for both smooth and non-smooth cases, to address this challenge.
The FPR method first establishes a homomorphism between the lower-dimensional definition domain of the quasiperiodic function and the higher-dimensional torus, then recovers the global quasiperiodic system by employing interpolation technique with finite points in the definition domain without dimensional lifting.
Furthermore, we develop accurate and efficient strategies of selecting finite points according to the arithmetic properties of irrational numbers. 
The corresponding mathematical theory, convergence analysis, and computational complexity analysis on choosing finite points are presented. 
Numerical experiments demonstrate the effectiveness and superiority of FPR approach in recovering both smooth quasiperiodic functions and piecewise constant Fibonacci quasicrystals. 
While existing spectral methods encounter difficulties in accurately recovering non-smooth quasiperiodic functions.
\end{abstract}



\section{Introduction}
Quasiperiodic systems, related to irrational numbers, have attracted extensive attention due to their fascinating mathematical properties~\cite{senechal1996quasicrystals, bohr2018almost, meyer2000algebraic, penrose1974role, steurer2009crystallography, baake2013aperiodic}. 
Quasiperiodic behavior is widely observed in physics and materials sciences, such as many-body celestial systems, quasicrystals, incommensurate systems, polycrystalline materials, and quantum systems~\cite{poincare1890probleme, shechtman1984metallic,cao2018unconventional, sutton1995interfaces, zeng2004supramolecular, hofstadter1976energy}. 
Among all quasiperiodic systems, the non-smooth case is of particular interest, such as discrete Schrödinger operator with quasiperiodic potential, Fibonacci photonic quasicrystal, and discrete time quasicrystal~\cite{avila2009ten, damanik1999uniform, merlin1985quasiperiodic, vardeny2013optics, tanese2014fractal, giergiel2019discrete, goblot2020emergence, verbin2015topological}.

Quasiperiodic systems pose significant challenges for numerical computation and corresponding theoretical analysis, due to their space-filling order without decay. 
Several numerical methods have been developed to address quasiperiodic systems. 
A widely used method, periodic approximation method (PAM)~\cite{zhang2008efficient}, employs periodic systems to approximate quasiperiodic systems over a finite domain, corresponding to  
using rational numbers to approximate irrational numbers in reciprocal space. 
PAM only approximates partial quasiperiodic systems and 
inevitably brings the rational approximation error, 
unless the period becomes infinity~\cite{jiang2023on}. 
Numerical examples have demonstrated that the rational approximation error plays a dominate role in numerically computing quasiperiodic systems~\cite{jiang2014numerical, jiang2015stability}.
An accurate algorithm is the quasiperiodic spectral method (QSM)~\cite{jiang2018numerical}.
Based on the Fourier-Bohr transformation, QSM can expand quasiperiodic functions with trigonometric polynomials. 
Theoretical analysis has been shown that QSM has exponential convergence for smooth cases~\cite{jiang2022numerical}. 
However, when dealing with nonlinear problems, the computational cost of QSM becomes unaffordable due to the unavailability of fast Fourier transform (FFT).
Another accurate approach is the projection method (PM) that captures the essential feature of quasiperiodic systems which can be embedded into higher-dimensional periodic systems~\cite{jiang2014numerical}. 
PM is an extension of the periodic Fourier pseudo-spectral method. 
The spectral allocation technique is employed to represent quasiperiodic functions by introducing the discrete Fourier-Bohr transformation. PM has exponential convergence for smooth systems and can utilize FFT to reduce computational complexity~\cite{jiang2022numerical}. 

\begin{figure*}[!hbpt]
    \centering
	\subfigure[Exact value]{\label{fig:fibonacci6_1}\includegraphics[width=4.5cm]{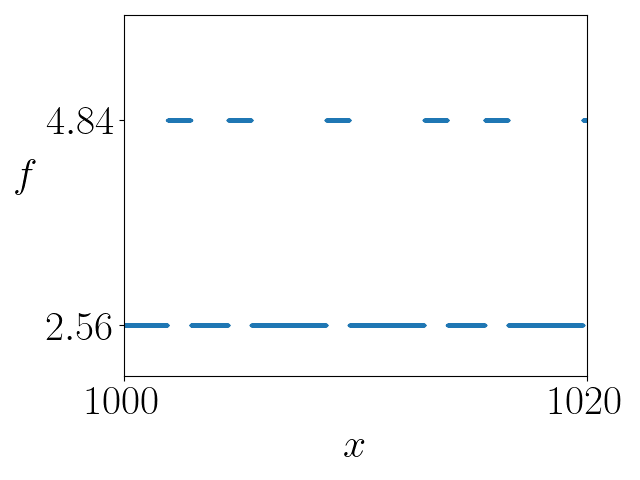}}
	\hspace{10mm}
	\subfigure[Result solved by PM method]{\label{fig:fibonacci6_3}\includegraphics[width=4.5cm]{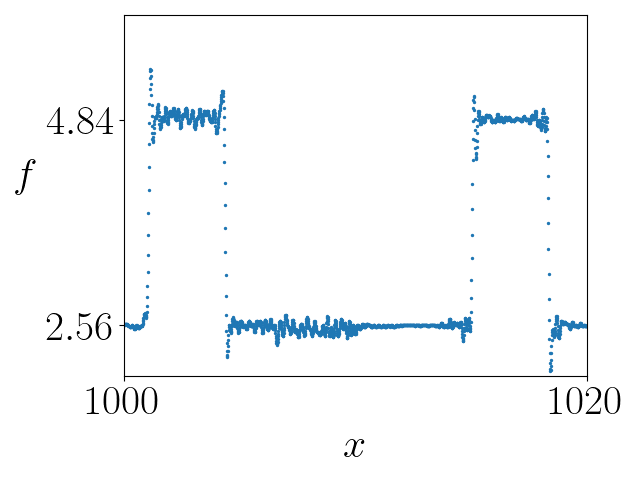}}
	\caption{Comparison of the exact value and the numerical result obtained by PM method for the piecewise dielectric function when $x\in[1000,1020)$.}
    \label{fig:fibonacci6}
\end{figure*}

Existing Galerkin spectral approaches, especially PM, have made progress in numerically solving quasiperiodic systems, including quasicrystals~\cite{jiang2014numerical, jiang2015stability}, incommensurate quantum systems~\cite{xueyang2021numerical,gao2023pythagoras}, topological insulators~\cite{wang2022effective}, grain boundaries~\cite{cao2021computing,jiang2022tilt}. 
However, none of these methods is suitable for solving non-smooth quasiperiodic problems. 
For example, consider the piecewise dielectric function of a 1D Fibonacci photonic quasicrystal (see \Cref{exa:fibonacci} for details). 
\Cref{fig:fibonacci6} shows the numerical result obtained by PM method, which is completely inconsistent with the exact value. 
Moreover, the Gibbs-phenomenon appears in a neighbourhood of the jump discontinuity. 
Hence, accurately recovering global non-smooth quasiperiodic systems remains an open problem, which motivates the development of new numerical algorithms. 

In this work, we pay attention to developing a new algorithm for recovering both smooth and non-smooth quasiperiodic systems. 
Our contributions are summarized as follows.
\begin{itemize}
    \item We propose a new approach, finite points recovery (FPR) method, for accurately recovering arbitrary dimensional quasiperiodic systems. 
    A homomorphism between the lower-dimensional definition domain of the quasiperiodic function and the higher-dimensional torus is established. 
    Based on this homomorphism, FPR method recovers the global quasiperiodic system by employing interpolation technique with finite points in the definition domain without dimensional lifting.
    \item We classify quasiperiodic systems into two categories according to the arithmetic properties of irrational numbers: badly approximable systems and good approximable systems. For each category, we employ distinct strategies for finite point selection within the FPR method, to ensure the accuracy and efficiency in the recovery process. 
    \item We provide a detailed exposition of the mathematical theory underlying the FPR method, along with rigorous proofs. Moreover, we present the convergence analysis of the algorithm and the computational complexity analysis on choosing finite points.
    \item We apply the FPR method to recover two classes of quasiperiodic systems, including smooth quasiperiodic functions and piecewise constant Fibonacci quasicrystals.
    Numerical experiments demonstrate the effectiveness and superiority of FPR approach in recovering the above two classes, while PM method fails to handle non-smooth systems.
\end{itemize}

The rest of this paper is organized as follows. 
In \Cref{sec:pre},  we give necessary notations and preliminary knowledge. 
In \Cref{sec_theo}, we establish a homomorphism between the lower-dimensional definition domain of the quasiperiodic function and the higher-dimensional torus.
In \Cref{sec:alg}, we propose the FPR method, analyze its convergence and computational complexity, and discuss the impact of the arithmetic properties of irrational numbers.
In \Cref{sec:num}, we show the accuracy and superiority of FPR method in handing both smooth and non-smooth quasiperiodic systems.
In \Cref{sec:con}, we carry out the summary and give an outlook of future work.

\section{Preliminaries}\label{sec:pre}
In this section, we introduce the requisite notations and preliminary knowledge. 
$\bm{I}_{n}$ is the $n$-order identity matrix. 
$\mathbb{N}_+$ denotes the set of all positive integers. 
The Cartesian product of two sets, $X$ and $Y$, denoted by $X\times Y$, is the set of all ordered pairs $(x,y)$, where $x$ and $y$ are elements of $X$ and $Y$, respectively. 
$X^d$ denotes the Cartesian product of $d$ sets $X$. 
$\mathbb{T}^n=\mathbb{R}^n/\mathbb{Z}^n$ is the $n$-dimensional torus. 
$\rm{rank}_{\mathbb{R}}$ ($\rm{rank}_{\mathbb{Q}}$) denotes the rank of a set of vectors or a matrix over the number field $\mathbb{R}$ ($\mathbb{Q}$).
$\dim_{\mathbb{R}}$ ($\dim_{\mathbb{Q}}$) denotes the dimension of a space over  $\mathbb{R}$ ($\mathbb{Q}$). 
For a vector $\bm{\alpha}=(\alpha_i)_{i=1}^n\in \mathbb{R}^n$, 
$[\bm{\alpha}]$ is the round down symbol of $\bm{\alpha}$ in each coordinate variable,  
and the infinity norm of $\bm{\alpha}$ is defined by $\|\bm{\alpha}\|_{\infty}:=\max_{1\leq i\leq n}|\alpha_i|$. 
For a region $\mathcal{H}=\{\bm{s}=(s_i)^n_{i=1}:\alpha_i\leq s_i\leq \beta_i\}\subset \mathbb{R}^n$, the vector $\langle\mathcal{H}\rangle:=(\beta_i-\alpha_i)_{i=1}^n$ measures the size of $\mathcal{H}$ and $\langle\mathcal{H}\rangle_i$ denotes its $i$-th component. 

For an $n$-dimensional periodic function $F(\bm{s})$, there exists a set of primitive vectors $\{\bm{a}_1,\cdots,\bm{a}_n\}$ that forms a basis of $\mathbb{R}^n$. The corresponding Bravais lattice is 
\begin{equation*}
    \mathcal{B}:=\{\bm{R}\in \mathbb{R}^n: \bm{R}=\ell_1\bm{a}_1+\cdots+\ell_n\bm{a}_n,~\bm{\ell}=(\ell_i)_{i=1}^n\in \mathbb{Z}^n\}.
\end{equation*}
For each $\bm{R}\in \mathcal{B}$, $F(\bm{s})$ satisfies $F(\bm{s}+\bm{R})=F(\bm{s})$. 
The unit cell of $F(\bm{s})$, denoted by $\Omega$, is the fundamental domain 
\begin{equation*}
    \Omega:=\{\bm{s}=w_1\bm{a}_1+\cdots+w_n\bm{a}_n: \bm{w}=(w_i)_{i=1}^n\in[0,1)^n\}.
\end{equation*}
There always exists an invertible linear transformation such that the basis $\{\bm{a}_1,\cdots,\bm{a}_n\}$ can be transform into a standard orthonormal basis $\{\bm{e}_1,\cdots,\bm{e}_n\}$. 
This enables us to give all the theoretical analysis in this paper based on the assumption that the unit cell is a cube $\Omega=[0,1)^n$.

Before we proceed, let's provide some required definitions.

\begin{definition}
    A matrix $\bm{P}\in \mathbb{R}^{d\times n}$ is the  projection matrix, if it belongs to the set $ \mathbb{P}^{d\times n}$ defined as
    \begin{equation*}
        \mathbb{P}^{d\times n}:=\{\bm{P}=(\bm{p}_1,\cdots,\bm{p}_n)\in\mathbb{R}^{d\times n}:\mathrm{rank}_{\mathbb{R}}(\bm{P})=d,~\mathrm{rank}_{\mathbb{Q}}(\bm{p}_1,\cdots,\bm{p}_n)=n\}.
    \end{equation*}
\end{definition}

\begin{definition}
    A $d$-dimensional function $f(\bm{x})$ is  quasiperiodic, 
    if there exists an $n$-dimensional periodic function $F(\bm{s})$ and a projection matrix $\bm{P}\in \mathbb{P}^{d\times n}$, such that $f(\bm{x})=F(\bm{P}^T\bm{x})$ for all $\bm{x}\in\mathbb{R}^d$. $F(\bm{s})$ is called the parent function of $f(\bm{x})$, and we refer to $\mathbb{R}^d$ as the physical space and $\mathbb{R}^n$ as the superspace. 
\end{definition}

\begin{definition}\label{def:lift_map}
    Let $\mathcal{G}$ be a set in the physical space $\mathbb{R}^d$ and $\bm{P}\in \mathbb{P}^{d\times n}$ be a projection matrix.
    The lift map $\mathcal{L}$ is defined by
    \begin{equation*}
        \begin{aligned}
        &\mathcal{L}: &\mathcal{G}&\to \mathbb{R}^n,\\
        &&\bm{x}&\mapsto \bm{P}^T\bm{x}.
        \end{aligned}
    \end{equation*}
\end{definition}

\begin{remark}\label{rem:slice}
    Let $\mathcal{S}:=\mathcal{L}(\mathbb{R}^d)$ be a subspace of $\mathbb{R}^n$. From the definition of the lift map $\mathcal{L}$, $\dim_{\mathbb{R}}(\mathcal{S})=d$ and $\dim_{\mathbb{Q}}(\mathcal{S})=n$. For convenience, we refer to $\mathcal{S}$ as an irrational slice in superspace $\mathbb{R}^n$. It is straightforward to show that $\mathcal{L}$ is an isomorphism from $\mathbb{R}^d$ to $\mathcal{S}$.
\end{remark}

\begin{remark}
    The irrational slice $\mathcal{S}$ contains one and only one lattice point $\bm{R}\in\mathcal{B}$ (see Chapter 2.5 in \cite{meyer2000algebraic}).  For convenience, we set $\bm{0}=(0,\cdots,0)\in\mathcal{S}$. For general cases, only one translation operation is required for the entire Bravais lattice. 
\end{remark}

\begin{definition}\label{def:mod_map}
    For any set $\mathcal{H}$ of superspace $\mathbb{R}^n$, 
    the modulo map $\mathcal{M}$ is defined by
    \begin{equation*}
        \begin{aligned}
        &\mathcal{M}: &\mathcal{H}&\to \mathbb{T}^n,\\
        &&\bm{s}&\mapsto \overline{\bm{s}},
        \end{aligned}
    \end{equation*}
    where $\overline{\bm{s}}:=\bm{s}+\mathbb{Z}^n$ is a left coset of $\mathbb{Z}^n$. We denote the image of $\mathcal{H}$ under $\mathcal{M}$ by $\overline{\mathcal{H}}:=\mathcal{M}(\mathcal{H})=\mathcal{H}/\mathbb{Z}^n$.
\end{definition}

\begin{remark}
    The modulo map $\mathcal{M}$ is a natural homomorphism from $\mathbb{R}^n$ to $\mathbb{T}^n$, 
    with $\overline{\mathbb{R}^n}=\mathbb{R}^n/\mathbb{Z}^n=\mathbb{T}^n$ and $\rm{Ker}(\mathcal{M})=\mathbb{Z}^n$. 
\end{remark}

\begin{remark}
    $\mathbb{T}^n$ and $\Omega$ are equivalent under the sense of isomorphism. 
    Specifically, we can define an isomorphism $\Phi:\mathbb{T}^n\rightarrow\Omega$ that maps each coset of $\mathbb{T}^n$ to its representative in $\Omega$. 
\end{remark}

\section{$\overline{\mathcal{S}}$ is dense in $\mathbb{T}^n$}\label{sec_theo}
In this section, we first observe that the $n$-dimensional irrational slice $\mathcal{S}$, after the modulo operation, is dense in $\mathbb{T}^n$. We then provide a rigorous proof of this observation. Next, we introduce a homomorphism between $\mathbb{R}^d$ and $\mathbb{T}^n$. Finally, we establish a close relationship between the arithmetic properties of irrational numbers and Diophantine approximation systems.

\subsection{Observation}\label{subsec:observation}

In this subsection, we concern the distribution of $\overline{\mathcal{S}}=\mathcal{M}(\mathcal{S})$ in $\mathbb{T}^n$.  
To describe the process of the modulo map $\mathcal{M}$ acting on $\mathcal{S}$, we present an equivalent expression of $\mathcal{S}$. 
Since $\dim_{\mathbb{R}}(\mathcal{S})=d$, we can decompose each point $\bm{s}\in \mathcal{S}$ as $\bm{s}=(\bm{t},\bm{r})^T$, 
where $\bm{t}=(t_i)_{i=1}^d$, $\bm{r}=(r_i)_{i=1}^{n-d}$,   
and $r_i~(i=1,\cdots,n-d)$ is a linear function of $t_1,\cdots,t_d$, i.e., 
\begin{equation*}
    r_i=r_i(t_1,\cdots,t_d)=\alpha_{i1}t_1+\cdots+\alpha_{id}t_d,
\end{equation*} 
where $\alpha_{ij}\in\mathbb{R},~j=1,\cdots,d$.
Let us introduce matrices $\bm{A} \in \mathbb{R}^{(n-d)\times d}$ and $\bm{Q} \in \mathbb{R}^{n\times d}$, where
\begin{equation*}\label{eq:matrix_a}
    \bm{A}=
    \begin{bmatrix}
        \alpha_{11}&\cdots &\alpha_{1d}\\
        \vdots& &\vdots\\
        \alpha_{n-d,1}&\cdots &\alpha_{n-d, d}\\
    \end{bmatrix},
\end{equation*}
and
\begin{equation*}\label{eq:matrix_q}
    \bm{Q}=
    \begin{bmatrix}
        \bm{I}_d \\ \bm{A}
    \end{bmatrix}. 
\end{equation*}
Using these two matrices, we can express the slice $\mathcal{S}$ as
\begin{equation}\label{eq:rewrite_slice}
	\begin{aligned}
		\mathcal{S}=&~\{\bm{s}=(\bm{t},\bm{r})^T : \bm{r}=\bm{A}\bm{t},~ \bm{t}\in \mathbb{R}^d\},\\
		=&~\{\bm{s}=(\bm{t},\bm{r})^T : \bm{s}=\bm{Q}\bm{t}, ~\bm{t}\in \mathbb{R}^d\}.\\
	\end{aligned}
\end{equation}
\begin{remark}
    The matrix $\bm{Q}$ has similar properties as $\bm{P}^T$. Specifically, its row vectors $\bm{q}_1,\cdots,$ $\bm{q}_n\in \mathbb{R}^d$ satisfy ${\rm rank}_{\mathbb{Q}}(\bm{q}_1,\cdots,\bm{q}_n)=n$ and ${\rm rank}_{\mathbb{R}}(\bm{Q})=d$. In fact, $\bm{P}^T$ can be transformed into $\bm{Q}$ by a linear transformation over $\mathbb{Q}$.
\end{remark}

\begin{figure*}[!hbpt]
    \centering
	\subfigure[$\mathcal{S}$ before modulo.]
	{\includegraphics[width=4.0cm]{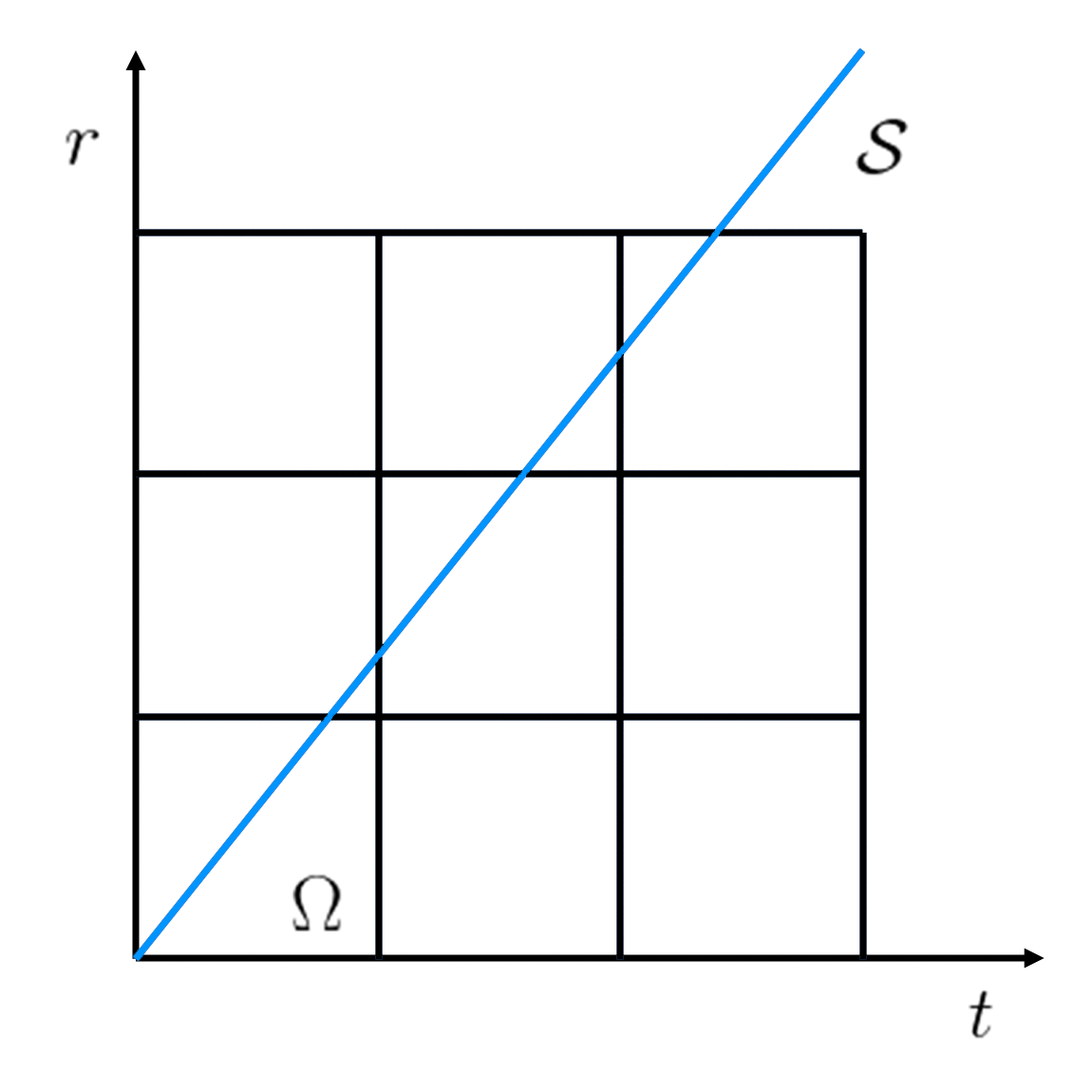}}
	\subfigure[Step 1: Modulo $\mathcal{S}$ along the $t$-axis.]
	{\label{fig:modulo(b)}\includegraphics[width=3.9cm]{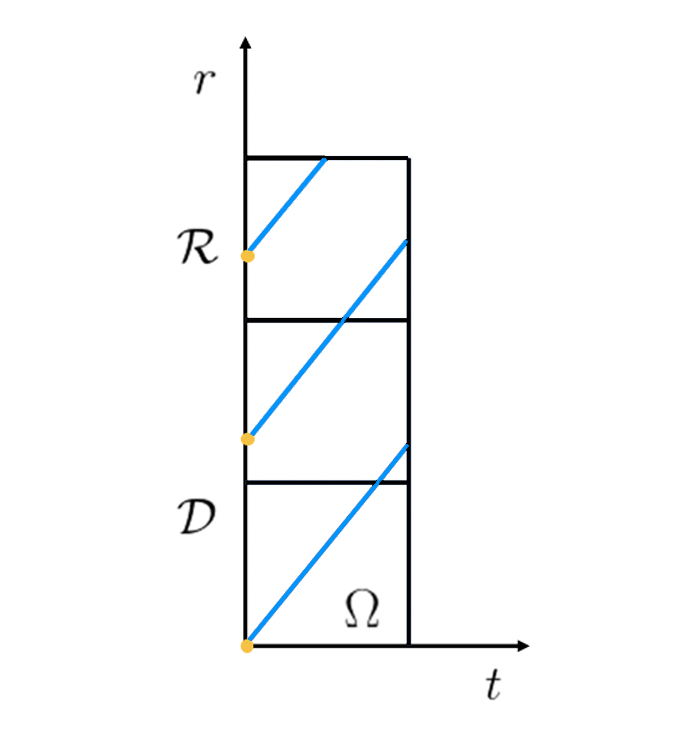}}
	\subfigure[Step 2: Modulo $\mathcal{S}$ along the $r$-axis after step 1.]
	{\label{fig:modulo(c)}\includegraphics[width=4.0cm]{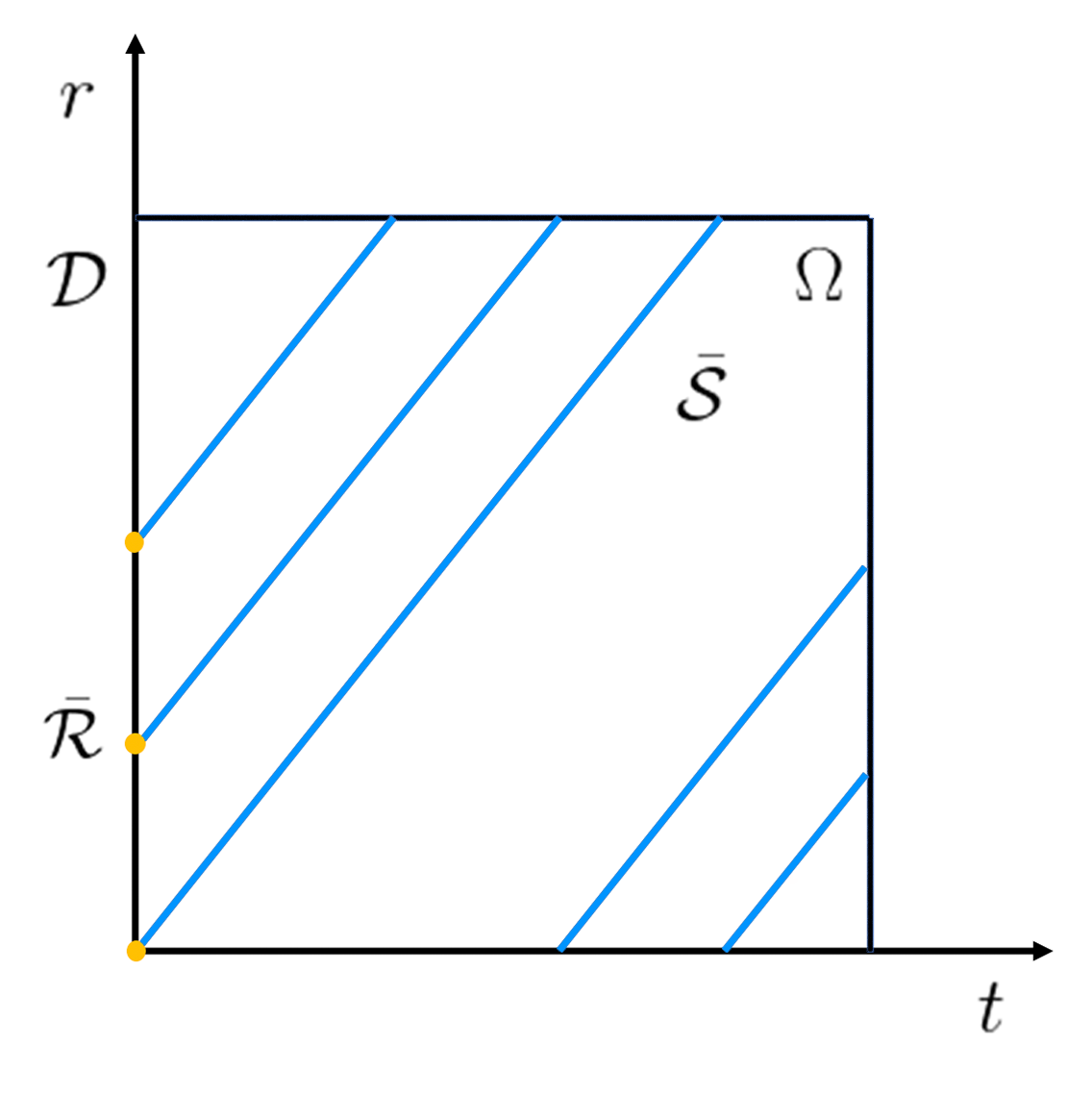}}	
	\caption{Illustration of the process of modulo a 2D irrational slice $\mathcal{S}=\{\bm{s}=(t,r)^T:r=\sqrt{2}t,t\in\mathbb{R}\}$. The images depict the modulo process within 9 unit cells.}
    \label{fig:modulo}
\end{figure*}

\begin{figure*}[!hbpt]
    \centering
	\subfigure[$t\in[0,10)$]{\includegraphics[width=4.0cm]{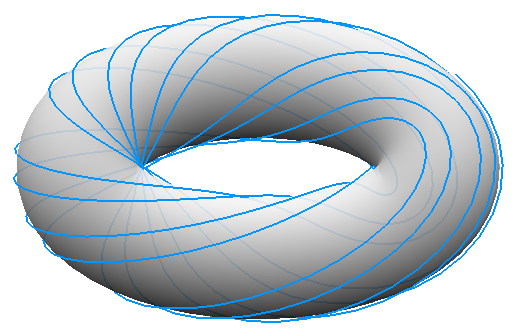}}
	\subfigure[$t\in[0,30)$]{\includegraphics[width=4.0cm]{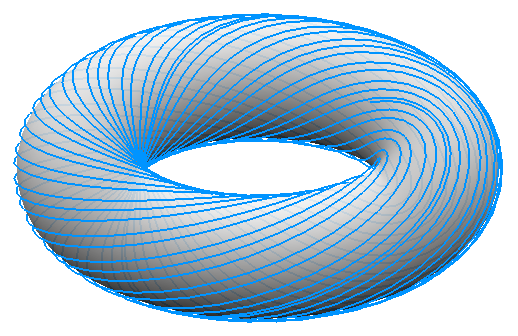}}
	\hspace{0.3cm}
	\subfigure[$t\in[0,70)$]{\includegraphics[width=4.0cm]{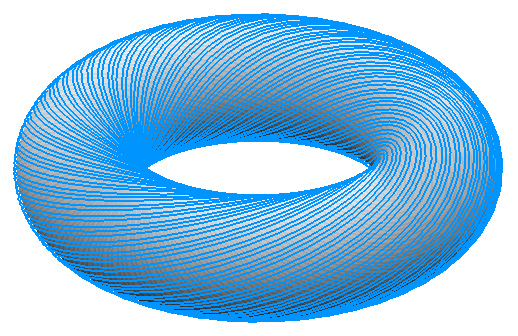}}	
	\caption{$\overline{\mathcal{S}}$ in $\mathbb{T}^2$, where irrational slice $\mathcal{S}=\{\bm{s}=(t,r)^T:r=\sqrt{2}t,t\in\mathbb{R}\}$.}
    \label{fig:torus}
\end{figure*}

We take an example to display the above observation. Consider a 2D irrational slice $\mathcal{S}=\{\bm{s}=(t,r)^T:r=\sqrt{2}t,~t\in\mathbb{R}\}$,
\Cref{fig:modulo} illustrates the process of modulo map $\mathcal{M}$ acting on $\mathcal{S}$ within 9 unit cells. 
The horizontal and vertical axes in the 2D plane are denoted as $t$-axis and $r$-axis, respectively. 
The modulo operation is divided into two steps, first along the $t$-axis (see \Cref{fig:modulo(b)}) and then along the $r$-axis (see \Cref{fig:modulo(c)}). 
The resulting slice family of $\mathcal{S}$ after the modulo operation is denoted as $\overline{\mathcal{S}}$. \Cref{fig:torus} shows that $\overline{\mathcal{S}}$ becomes denser in $\mathbb{T}^2$ as the range of $t$ increases.

Let $\mathcal{R}:=\{\bm{0}\}\times \{\bm{r}=\bm{A}\bm{t}:\bm{t}\in \mathbb{Z}^d\}$ 
and $\mathcal{D}:=\{\bm{0}\}\times\mathbb{T}^{n-d}$~($\bm{0}=(0,\cdots,0)^T\in \mathbb{R}^{d}$), 
then $\overline{\mathcal{R}}=\mathcal{R}/\mathbb{Z}^n
    =\{\bm{0}\}\times\{\bm{r}=\bm{A}\bm{t}:\bm{t}\in \mathbb{Z}^d\}/\mathbb{Z}^{n-d}=\overline{\mathcal{S}}\cap\mathcal{D}$. 
\Cref{fig:modulo(c)} shows that the distribution of $\overline{\mathcal{S}}$ in $\mathbb{T}^n$ can be completely determined by the distribution of $\overline{\mathcal{R}}$ in $\mathcal{D}$, i.e., the distribution of $\{\bm{r}=\bm{A}\bm{t}:\bm{t}\in \mathbb{Z}^d\}/\mathbb{Z}^{n-d}$ in $\mathbb{T}^{n-d}$. 
It inspires us to develop a rigorous mathematical theory of the observation ``$\overline{\mathcal{S}}$ is dense in $\mathbb{T}^n$", which is presented in the next subsection.

\subsection{Theoretical analysis}\label{subsec:math_theo}
In this subsection, we abstract the above observation into a mathematical theorem and present a rigorous proof. 
Then we introduce a homomorphism between $\mathbb{R}^d$ and $\mathbb{T}^n$. 
Let us first introduce this theorem.

\begin{theorem}\label{thm:dens_fill}
    For any irrational slice $\mathcal{S}\subset\mathbb{R}^n$,  
    $\overline{\mathcal{S}}=\mathcal{M}(\mathcal{S})$ is dense in $\mathbb{T}^n$.
\end{theorem}

Before the proof, we need to give the definition of dense subset in the space $\mathbb{R}^n$, as well as the $k$-variable Kronecker theorem.

\begin{definition}\label{def:dense}
    Let $V_1\subset V_2\subseteq \mathbb{R}^n$, the subset $V_1$ is dense in $V_2$ if for any $\bm{s}_2\in V_2$ and any $\epsilon>0$, there exists a $\bm{s}_1\in V_1$ such that $\|\bm{s}_1-\bm{s}_2\|_{\infty}<\epsilon$. 
\end{definition}

\begin{lemma}[$k$-variable Kronecker theorem~\cite{granville2020number}]\label{lem:kron}
    Given $k$ real numbers $\alpha_1,\cdots,$ $\alpha_k$, 
    assume that $1,\alpha_1,\cdots,\alpha_k$ are $\mathbb{Q}$-linearly 
    independent, then the point set
    $\{(\alpha_1,\cdots,\alpha_k)^Tm: m\in \mathbb{Z}\}/\mathbb{Z}^k$ is dense in $\mathbb{T}^k$.
\end{lemma}

With these preparations, we start the proof of \cref{thm:dens_fill}.

\begin{proof}
    First, let us establish the equivalence between ``$\overline{\mathcal{S}}=\mathcal{M}(\mathcal{S})$ is dense in $\mathbb{T}^n$" and ``$\{\bm{r}=\bm{A}\bm{t}:\bm{t}\in \mathbb{Z}^d\}/\mathbb{Z}^{n-d}$ is dense in $\mathbb{T}^{n-d}$". 
    From \Cref{def:dense}, ``$\overline{\mathcal{S}}=\mathcal{M}(\mathcal{S})$ is dense in $\mathbb{T}^n$" means that for each $\bm{s}^*=(\bm{t}^*,\bm{r}^*)^T\in\mathbb{T}^n$ and $\epsilon>0$, 
    there exists a $\bm{s}=(\bm{t},\bm{A}\bm{t})^T\in\mathcal{S}$ such that $\|\overline{\bm{s}}-\bm{s}^*\|_{\infty}<\epsilon$. 
    Let $\bm{b}^*:=\bm{r}^*-\bm{A}\bm{t}^*$, then $\bm{b}^*\in \mathbb{T}^{n-d}$. 
    The arbitrariness of $\bm{s}^*$ in $\mathbb{T}^n$ implies the arbitrariness of $\bm{b}^*$ in $\mathbb{T}^{n-d}$.
    Meanwhile, let $\bm{t}':=[\bm{t}]\in\mathbb{Z}^d$, then  $\bm{s}'=(\bm{t}',\bm{A}\bm{t}')^T\in\mathcal{S}$ satisfies 
    $\overline{\bm{A}\bm{t}'}=\overline{\bm{A}\bm{t}'}-\bm{A}\overline{\bm{t}'}=\overline{\bm{A}\bm{t}}-\bm{A}\overline{\bm{t}}:=\bm{b}$. 
    Note that, $\|\overline{\bm{s}}-\bm{s}^*\|_{\infty}<\epsilon$ implies $\|\bm{b}-\bm{b}^*\|_{\infty}<\epsilon$. 
    Then for each $\bm{b}^*\in \mathbb{T}^{n-d}$ and $\epsilon>0$, there exists a $\bm{t}'\in\mathbb{Z}^d$ such that $\|\overline{\bm{A}\bm{t}'}-\bm{b}^*\|_{\infty}<\epsilon$, which exactly means ``$\{\bm{r}=\bm{A}\bm{t}:\bm{t}\in \mathbb{Z}^d\}/\mathbb{Z}^{n-d}$ is dense in $\mathbb{T}^{n-d}$". Here, we have completed the proof for one direction of the equivalence, and the proof for the other direction can be obtained similarly from the above process. 
    
    Next, we use the induction method on the dimension $d$ of $\bm{t}$ to prove that ``$\{\bm{r}=\bm{A}\bm{t}:\bm{t}\in \mathbb{Z}^d\}/\mathbb{Z}^{n-d}$ is dense in $\mathbb{T}^{n-d}$". 
    When $d=1$, $\bm{A}=(\alpha_{11},\cdots,\alpha_{n-1,1})^T$, 
    $\bm{Q}=(1,\alpha_{11},\cdots,\alpha_{n-1,1})^T$, 
    and the elements of $\bm{Q}$ are $\mathbb{Q}$-linearly independent.
    By \Cref{lem:kron}, $\{(\alpha_{11},\cdots,\alpha_{n-1,1})^T m:m\in \mathbb{Z}\}/\mathbb{Z}^{n-1}$ 
    is dense in $\mathbb{T}^{n-1}$. 
    Hence, $\{\bm{r}=\bm{A}t:t\in \mathbb{Z}\}/\mathbb{Z}^{n-1}$ is dense in  $\mathbb{T}^{n-1}$.
    Assume that the conclusion holds when $d=k-1$, 
    we prove that it also holds for $d=k$ in two cases. 
    
    Case 1: Suppose that the $k$-th column of matrix $\bm{A}$ consists entirely of rational numbers, i.e., $\alpha_{1k},\cdots,\alpha_{n-k,k}\in\mathbb{Q}$. 
    	Note that 
    	\begin{equation*}
    		\begin{aligned}
    			\bm{r}=\bm{A}\bm{t}&=
    			\begin{bmatrix}
    				\alpha_{11} & \cdots & \alpha_{1,k-1} & \alpha_{1k}\\
    				\vdots & & \vdots & \vdots\\
    				\alpha_{n-k,1} & \cdots & \alpha_{n-k,k-1} & \alpha_{n-k,k}\\
    			\end{bmatrix} 
    			\begin{bmatrix}
    				t_1 \\ \vdots \\ t_k
    			\end{bmatrix} \\
    			&=\begin{bmatrix}
    				\alpha_{11} & \cdots & \alpha_{1,k-1} \\
    				\vdots & & \vdots \\
    				\alpha_{n-k,1} & \cdots & \alpha_{n-k,k-1} \\
    			\end{bmatrix} 
    			\begin{bmatrix}
    				t_1 \\ \vdots \\ t_{k-1}
    			\end{bmatrix}
    			+\begin{bmatrix}
    				\alpha_{1k} \\ \vdots \\ \alpha_{n-k,k}
    			\end{bmatrix}t_k. \\
    		\end{aligned}
    	\end{equation*}
    	Define matrices $\bm{A}_1 \in \mathbb{R}^{(n-k)\times(k-1)}$ and $\bm{Q}_1 \in \mathbb{R}^{(n-1)\times(k-1)}$ as follows
    	\begin{equation*}
    		\bm{A}_1:=
    		\begin{bmatrix}
    			\alpha_{11} & \cdots & \alpha_{1,k-1} \\
    			\vdots & & \vdots \\
    			\alpha_{n-k,1} & \cdots & \alpha_{n-k,k-1} \\
    		\end{bmatrix},\quad 
    		\bm{Q}_1:=\begin{bmatrix}
    			\bm{I}_{k-1} \\ \bm{A}_1
    		\end{bmatrix}.
    	\end{equation*}
    	Let $\bm{t}_1=(t_i)_{i=1}^{k-1}$ and $\bm{r}_1=\bm{A}_1\bm{t}_1$.
    	Since $\alpha_{1k},\cdots,\alpha_{n-k,k}$ are rational, 
    	$1$ and $\alpha_{i k}$ are $\mathbb{Q}$-linearly dependent for all $i\in \{1,\cdots,n-k\}$. 
    	This implies that the row vectors of $\bm{Q}_1$ are $\mathbb{Q}$-linearly independent. 
    	Otherwise, suppose that there exists $i\in \{1,\cdots,n-k\}$ 
    	such that the $i$-th row of $\bm{A}_1$ and all rows of $\bm{I}_{k-1}$ are $\mathbb{Q}$-linearly dependent. 
    	Meanwhile, $1$ and $\alpha_{i k}$ are also $\mathbb{Q}$-linearly dependent. It follows that the $i$-th row of $\bm{A}$ and all rows of $\bm{I}_{k}$ are $\mathbb{Q}$-linearly dependent, which is contradictory to the properties of $\bm{Q}$.
    	According to the inductive hypothesis, 
    	$\{\bm{r}_1=\bm{A}_1\bm{t}_1: \bm{t}_1\in \mathbb{Z}^{k-1}\}/\mathbb{Z}^{n-k}$
    	is dense in $\mathbb{T}^{n-k}$, 
    	then $\{\bm{r}=\bm{r}_1+(\alpha_{1k},\cdots,\alpha_{n-k,k})^T t_k: t_k\in \mathbb{Z}\}/\mathbb{Z}^{n-k}$ is dense in $\mathbb{T}^{n-k}$.
    
    Case 2: Suppose that in the $k$-th column of $\bm{A}$, there are at most $q$ elements that are $\mathbb{Q}$-linearly independent with $1$. Without loss of generality, we assume that $\alpha_{n-k-q+1,k},\cdots,\alpha_{n-k,k}$ are the $q$ linearly independent elements in the $k$-th column of $\bm{A}$. Denote
    	\begin{equation*}
    		\bm{A}_1:=
    		\begin{bmatrix}
    			\alpha_{11} & \cdots & \alpha_{1,k-1} & 0\\
    			\vdots & & \vdots \\
    			\alpha_{n-k-q,1} & \cdots & \alpha_{n-k-q,k-1} & 0\\
    			0 & \cdots & 0 & \alpha_{n-k-q+1,k}\\
    			\vdots & & \vdots & \vdots \\
    			0 & \cdots & 0 & \alpha_{n-k,k}
    		\end{bmatrix}\in \mathbb{R}^{(n-k)\times k} ,
    	\end{equation*}
    	and 
    	\begin{equation*}
    		\bm{A}_2:=
    		\begin{bmatrix}
    			\alpha_{n-k-q+1,1} & \cdots & \alpha_{n-k-q+1,k-1}\\
    			\vdots & & \vdots\\
    			\alpha_{n-k,1} & \cdots & \alpha_{n-k,k-1}
    		\end{bmatrix}\in \mathbb{R}^{q\times (k-1)}.
    	\end{equation*}
    	Then, we can express $\bm{r}$ as 
    	\begin{equation*}
    		\bm{r}=\bm{A}\bm{t}=\bm{A}_1\bm{t}+\bm{A}_2
    		\begin{bmatrix}
    			t_{1} \\ \vdots \\ t_{k-1}
    		\end{bmatrix}
    		+\begin{bmatrix}
    			\alpha_{1k} \\ \vdots \\ \alpha_{n-k-q,k}
    		\end{bmatrix}t_{k}.
    	\end{equation*}
    	Let $\bm{r}_1=(\bm{r}_{11},\bm{r}_{12})^T:=\bm{A}_1\bm{t}$, where
    	\begin{equation*}
    		\bm{r}_{11}
    		=\begin{bmatrix}
    			\alpha_{11} & \cdots & \alpha_{1,k-1} \\
    			\vdots & & \vdots \\
    			\alpha_{n-k-q,1} & \cdots & \alpha_{n-k-q,k-1} \\
    		\end{bmatrix}
    		\begin{bmatrix}
    			t_1 \\ \vdots \\ t_{k-1}
    		\end{bmatrix}:=\bm{A}_{11}\bm{t}_1,
    	\end{equation*}
    	\begin{equation*}
    		\bm{r}_{12}
    		=\begin{bmatrix}
    			\alpha_{n-k-q+1,k} \\ \vdots \\ \alpha_{n-k,k}
    		\end{bmatrix}t_k.
    	\end{equation*}
    	Let 
    	\begin{equation*}
    	    \bm{Q}_{11}:=\begin{bmatrix}
    	        \bm{I}_{k-1} \\ \bm{A}_{11}
    	    \end{bmatrix}\in \mathbb{R}^{(n-q-1)\times (k-1)},
    	\end{equation*}
    	then the row vectors of $\bm{Q}_{11}$ are $\mathbb{Q}$-linearly independent. 
    	Otherwise, the same reason as in Case 1 gives a contradiction. 
    	Again from the inductive hypothesis, 
    	$\{\bm{r}_{11}=\bm{A}_{11}\bm{t}_1: \bm{t}_1\in \mathbb{Z}^{k-1}\}/\mathbb{Z}^{n-k-q}$ is dense in $\mathbb{T}^{n-k-q}$.
    	Moreover, $1,\alpha_{n-k-q+1,k},\cdots,\alpha_{n-k,k}$ 
    	are $\mathbb{Q}$-linearly independent. 
    	From \Cref{lem:kron},  
    	$\{\bm{r}_{12}=(\alpha_{n-k-q+1,k},\cdots,\alpha_{n-k,k})^Tt_{k}: ~t_{k}\in \mathbb{Z}\}/\mathbb{Z}^q$ is dense in $\mathbb{T}^{q}$.
    	Therefore,  $\{\bm{r}_1=(\bm{r}_{11},\bm{r}_{12})^T\}/\mathbb{Z}^{n-k}$ is dense in $\mathbb{T}^{n-k}$, 
    	and then 
    	$\{\bm{r}=\bm{r}_1+\bm{A}_2\bm{t}_1+(\alpha_{1k},\cdots,\alpha_{n-k-q,k})^T t_k\}/\mathbb{Z}^{n-k}$
    	is dense in $\mathbb{T}^{n-k}$.
    	
    To summarize, we have proven that $\{\bm{r}=\bm{A}\bm{t}: \bm{t}\in \mathbb{Z}^d\}/\mathbb{Z}^{n-d}$ is dense in $\mathbb{T}^{n-d}$. Consequently, $\overline{\mathcal{S}}$ is dense in $\mathbb{T}^n$. 
%   Thus, we have completed the proof.
\end{proof}

Further, we introduce a homomorphism between $\mathbb{R}^d$ and $\mathbb{T}^n$.

\begin{definition}\label{def:comb_map}
    A combination map $\mathcal{C}$ is defined as the composition of the lift map $\mathcal{L}:\mathcal{G}\to \mathbb{R}^n$ and the modulo map $\mathcal{M}:\mathbb{R}^n\to \mathbb{T}^n$, denoted by
    \begin{equation*}
        \mathcal{C}:= \mathcal{M} \circ \mathcal{L},
    \end{equation*}
    where the symbol $\circ$ represents the composition of two maps.
\end{definition}

\begin{remark}
The combination map $\mathcal{C}$ is a homomorphism from $\mathbb{R}^d$ to $\mathbb{T}^n$. 
Moreover, from \Cref{thm:dens_fill}, $\mathcal{C}(\mathbb{R}^d)$ is dense in $\mathbb{T}^n$.
\end{remark}

Based on this homomorphism, we can obtain the following theorem which plays a crucial role in our proposed FPR method. 

\begin{theorem}\label{thm:relation}
    If $f(\bm{x})$ is a $d$-dimensional quasiperiodic function, there exists a parent function $F(\bm{s})$ and a projection matrix $\bm{P}\in \mathbb{P}^{d\times n}$, such that $f(\bm{x})=F(\bm{P}^T\bm{x})$ for all $\bm{x}\in\mathbb{R}^d$.
    Then 
    \begin{equation*}
        f(\bm{x})=F\left(\mathcal{C}(\bm{x})\right),\quad \forall\bm{x}\in \mathbb{R}^d,
    \end{equation*}
    where the combination map $\mathcal{C}$ is defined by \Cref{def:comb_map}.
\end{theorem}
\begin{proof}
    From the definition of $\mathcal{L}$, 
    \begin{equation*}
        f(\bm{x})=F(\mathcal{L}(\bm{x})),\quad \forall\bm{x}\in \mathbb{R}^d.
    \end{equation*}
    Further, from the periodicity of $F(\bm{s})$, 
    \begin{equation*}
        F(\bm{s})=F(\mathcal{M}(\bm{s})),\quad \forall\bm{s}\in \mathbb{R}^n.
    \end{equation*}
    Therefore, 
    \begin{equation*}
        f(\bm{x})=F(\mathcal{M}(\mathcal{L}(\bm{x})))=F(\mathcal{C}(\bm{x})),\quad \forall\bm{x}\in \mathbb{R}^d.
    \end{equation*}
\end{proof}

Using the homomorphism $\mathcal{C}$ and \Cref{thm:relation}, we can locate the point in a unit cell of parent function, corresponding to a point in physical space, and obtain the value of parent function at it.
This allows us to directly carry out our 
algorithm in the physical space based on finite points, without the need of dimensional lifting.

\subsection{ $\overline{\mathcal{S}}=\mathbb{T}^{n}$ or $\overline{\mathcal{S}}\neq \mathbb{T}^{n}$\label{subsec:bad_appr}}
In this subsection, we further explore whether $\overline{\mathcal{S}}=\mathbb{T}^{n}$ or not. 
From the proof of \Cref{thm:dens_fill}, this question is equivalent to whether $\{\bm{r}=\bm{A}\bm{t}:\bm{t}\in \mathbb{Z}^d\}/\mathbb{Z}^{n-d}=\mathbb{T}^{n-d}$ holds true or not. It is closely related to the Diophantine approximation problem.  
For a given matrix $\bm{A}\in\mathbb{R}^{k\times d}$, and for each $\varepsilon>0$, $\bm{r}_0\in \mathbb{T}^k$, 
we want to know whether there exists a solution $\bm{t}\in\mathbb{Z}^d$ satisfying the $k$-dimensional \textit{Diophantine system}
\begin{equation}\label{eq:diop_ineq}
    \|\bm{A}\bm{t}/\mathbb{Z}^k-\bm{r}_0\|_{\infty}\leq \varepsilon.
\end{equation}

\begin{remark}
    For a $d$-dimentional quasiperiodic function $f(\bm{x})$, if its projection matrix $\bm{P}\in\mathbb{R}^{d\times n}$ fulfils that $\bm{P}^T$ can be linearly transformed over $\mathbb{Q}$ into the form $\bm{Q}=(\bm{I}_d,\bm{A})^T$, then we refer to $f(\bm{x})$ as be related to the $(n-d)$-dimensional Diophantine system \cref{eq:diop_ineq}.
\end{remark}

\begin{definition}
    For each positive $\varepsilon$, $\mathcal{G}(\varepsilon)$ is the \textit{least region} of the Diophantine system \cref{eq:diop_ineq}, if $\langle\mathcal{G}(\varepsilon)\rangle\in\mathbb{N}^d_+$ and $\mathcal{G}(\varepsilon)$ contains a solution $\bm{t}\in \mathbb{Z}^d$ of \cref{eq:diop_ineq}. 
\end{definition}

\begin{lemma}(Chapter 1.3 in \cite{meyer2000algebraic})
\label{lem:lower_bound}
    For each matrix $\bm{A}\in\mathbb{R}^{k\times d}$, the least region $\mathcal{G}(\varepsilon)$ satisfies
    \begin{equation*}
        \langle\mathcal{G}(\varepsilon)\rangle_i\geq \frac{1}{2\varepsilon^k},~~i=1,\cdots,d.
    \end{equation*}
\end{lemma}

Next, we discuss the Diophantine approximation systems from two aspects, badly approximable systems and good approximable systems, based on the arithmetic properties of irrational numbers in $\bm{A}$.

\subsubsection{Badly approximable systems}

\begin{definition}\label{def:bad_appro}
    A matrix $\bm{A}\in\mathbb{R}^{k\times d}$ is \textit{badly approximable}, if there exists a sequence of positive integers $C_1,\cdots,C_d$, such that the least region $\mathcal{G}(\varepsilon)$ fulfills
    \begin{equation*}
        \frac{1}{2\varepsilon^{k}}\leq\langle\mathcal{G}(\varepsilon)\rangle_i\leq \frac{C_i}{\varepsilon^{k}},~~i=1,\cdots,d.
    \end{equation*}
    When $\bm{A}$ is a badly approximable matrix, we call \cref{eq:diop_ineq} a badly approximable system.
\end{definition}

\begin{remark}\label{rem:sqrt_2}
    When $k=1$, $\bm{A}$ becomes a badly approximable irrational number. 
    The set of badly approximable irrational numbers has the same cardinality as $\mathbb{R}$. In particular, all quadratic irrational numbers are badly approximable~\cite{burger2000exploring}.
\end{remark}

\begin{remark}\label{rem:rapid_fill}
    The upper bound of $\langle\mathcal{G}(\varepsilon)\rangle$ well controls the existence range of the solution, which is also known as the \textit{rapid filling} property of badly approximable systems~\cite{meyer2000algebraic}. 
\end{remark}

\begin{lemma}(Chapter 1.3 in \cite{meyer2000algebraic})
    When $\bm{A}$ is a badly approximable matrix, there exists a positive constant $C$ such that, for each $\bm{t}\in \mathbb{Z}^d\backslash\{\bm{0}\}$ and $\bm{r}_0\in\mathbb{T}^k$,
    \begin{equation}\label{eq:badly}
        \|\bm{A}\bm{t}/\mathbb{Z}^k-\bm{r}_0\|_{\infty}\geq \frac{C}{\|\bm{t}\|^k_{\infty}}>0.
    \end{equation}
\end{lemma}

For badly approximable systems, \cref{eq:badly} implies that there is no $\bm{t}\in\mathbb{Z}^d$ such that $\bm{A}\bm{t}/\mathbb{Z}^k=\bm{r}_0$ holds for any $\bm{r}_0\in \mathbb{T}^k$.
Therefore, when $k=n-d$, the first conclusion of this subsection comes out naturally.

\begin{theorem}
    When $\bm{A}\in\mathbb{R}^{(n-d)\times d}$ is a badly approximable matrix, 
    $\{\bm{r}=\bm{A}\bm{t}:\bm{t}\in \mathbb{Z}^d\}/\mathbb{Z}^{n-d}\neq\mathbb{T}^{n-d}$, i.e., $\overline{\mathcal{S}}\neq\mathbb{T}^n$.
\end{theorem}

\subsubsection{Good approximable systems}

\begin{definition}
    An irrational number $\alpha$ is \textit{good approximable}, if its continued fraction expansion $\alpha=[a_0,a_1,\cdots,a_n,\cdots]$ satisfies %$\{a_n\}_{n=1}^{\infty}$ is unbounded, i.e., 
    $$\overline{\lim_{n\to\infty}}a_n=\infty.$$
\end{definition}

\begin{remark}
    The set of good approximable numbers contains all Liouville numbers (see~\cite{liouville1844classes}
     for the definition of Liouville numbers). 
    Meanwhile, some famous transcendental numbers, like $e$ and $\pi$, are not Liouville numbers, but good approximation numbers. 
    Good approximable numbers form the complement of the set of badly approximable numbers within the set of irrational numbers~\cite{meyer2000algebraic}.
\end{remark}

\begin{definition}\label{def:good_appro}
    A matrix $\bm{A}\in\mathbb{R}^{k\times d}$ is \textit{good approximable}, 
    if all irrational numbers in matrix $\bm{A}$ are good approximable. 
    When $\bm{A}$ is good approximable, we call \cref{eq:diop_ineq} a good approximable system.
\end{definition}

Different from badly approximable systems, good approximable systems have arbitrary approximation property~\cite{hardy1979introduction}. 
It means that when $\bm{A}$ is good approximable, for each $\bm{r}_0\in\mathbb{T}^k$, there exists a $\bm{t}\in\mathbb{Z}^d$ such that $\bm{A}\bm{t}/\mathbb{Z}^k=\bm{r}_0$. Then we obtain the following conclusion when $k=n-d$.
\begin{theorem}
    When $\bm{A}\in\mathbb{R}^{(n-d)\times d}$ is a good approximable matrix, 
    $\{\bm{r}=\bm{A}\bm{t}:\bm{t}\in \mathbb{Z}^d\}/\mathbb{Z}^{n-d}=\mathbb{T}^{n-d}$, i.e., $\overline{\mathcal{S}}=\mathbb{T}^n$.
\end{theorem}

\section{Algorithm and analysis}\label{sec:alg}

From \Cref{thm:relation}, the global information of a quasiperiodic function is contained in the unit cell of parent function.
However, in practice, performing calculations in superspace may lead to an unbearable computational complexity.
In this section, we propose a new algorithm, FPR method, to recover the global quasiperiodic function by employing interpolation technique with finite points in physical space.  
We also present the convergence analysis of FPR method and the computational complexity of choosing finite points. 

\subsection{FPR method}

In this subsection, we present the FPR method in a step-by-step way. Assume that the projection matrix $\bm{P}\in\mathbb{P}^{d\times n}$ of the quasiperiodic function $f(\bm{x})$ is known, and $\bm{P}^T$ can be linearly transformed over $\mathbb{Q}$ into the form $\bm{Q}=(\bm{I}_d,\bm{A})^T$. 
Given a target point $\bm{x}^*\in\mathbb{R}^d$, FPR method has three steps to obtain an approximate value of $f(\bm{x}^*)$.

\textbf{Select region $\mathcal{G}$.}
FPR method uses finite points in physical space to recover the global quasiperiodic function.  
These finite points fall in a finite region $\mathcal{G}\subset\mathbb{R}^d$, which need to be determined in the first step. 
The selection principle of $\mathcal{G}$ is that $\mathcal{C}(\mathcal{G})$ is relatively uniform distributed in $\Omega$. 
Actually, the arithmetic properties of irrational numbers in the projection matrix $\bm{P}$ can provide a fast way for the selection of $\mathcal{G}$. We will give a more deep discussion in \Cref{subsec:rapid_filling}. 

\textbf{Shrink region $\mathcal{G}$.}
For the target point $\bm{x}^*$, 
its image under the lift map $\mathcal{L}$ is denoted as $\bm{s}^*=(\bm{t}^*, \bm{r}^*)^T$.
Define $\bm{b}^*:=\overline{\bm{A}[\bm{t}^*]}\in[0,1)^{n-d}$.  
Given $h>0$ and $\varepsilon>0$, we traverse the points in $\mathcal{G}\cap\mathbb{Z}^d$ to
find $n-d$ pairs of $(\bm{x}^{i+},\bm{x}^{i-})$, $i=1,\cdots,n-d$, such that
\begin{equation*}
    \begin{aligned}
    &|b^{i+}_i-(b^*_i+h/2)|\leq \varepsilon/2 &\text{and}\quad &|b^{i+}_j-b^*_i|\leq\varepsilon/2,\forall j\in\{1,\cdots,n-d\}\backslash \{i\},\\
    &|b^{i-}_i-(b^*_i-h/2)|\leq \varepsilon/2 &\text{and}\quad&|b^{i-}_j-b^*_i|\leq\varepsilon/2,\forall j\in\{1,\cdots,n-d\}\backslash \{i\},
    \end{aligned}
\end{equation*}
where 
\begin{equation*}
    \begin{aligned}
    &\bm{b}^{i+}:=\overline{\bm{A}\bm{t}^{i+}},\quad \bm{s}^{i+}=(\bm{t}^{i+},\bm{r}^{i+})^T:=\mathcal{L}(\bm{x}^{i+}),\\
    &\bm{b}^{i-}:=\overline{\bm{A}\bm{t}^{i-}},\quad \bm{s}^{i-}=(\bm{t}^{i-},\bm{r}^{i-})^T:=\mathcal{L}(\bm{x}^{i-}).\\
    \end{aligned}
\end{equation*}
Note that the traverse process in $\mathcal{G}\cap\mathbb{Z}^d$ can be performed simultaneously in each dimension.

\begin{figure*}[!hbpt]
    \centering
	\includegraphics[width=0.4\textwidth]{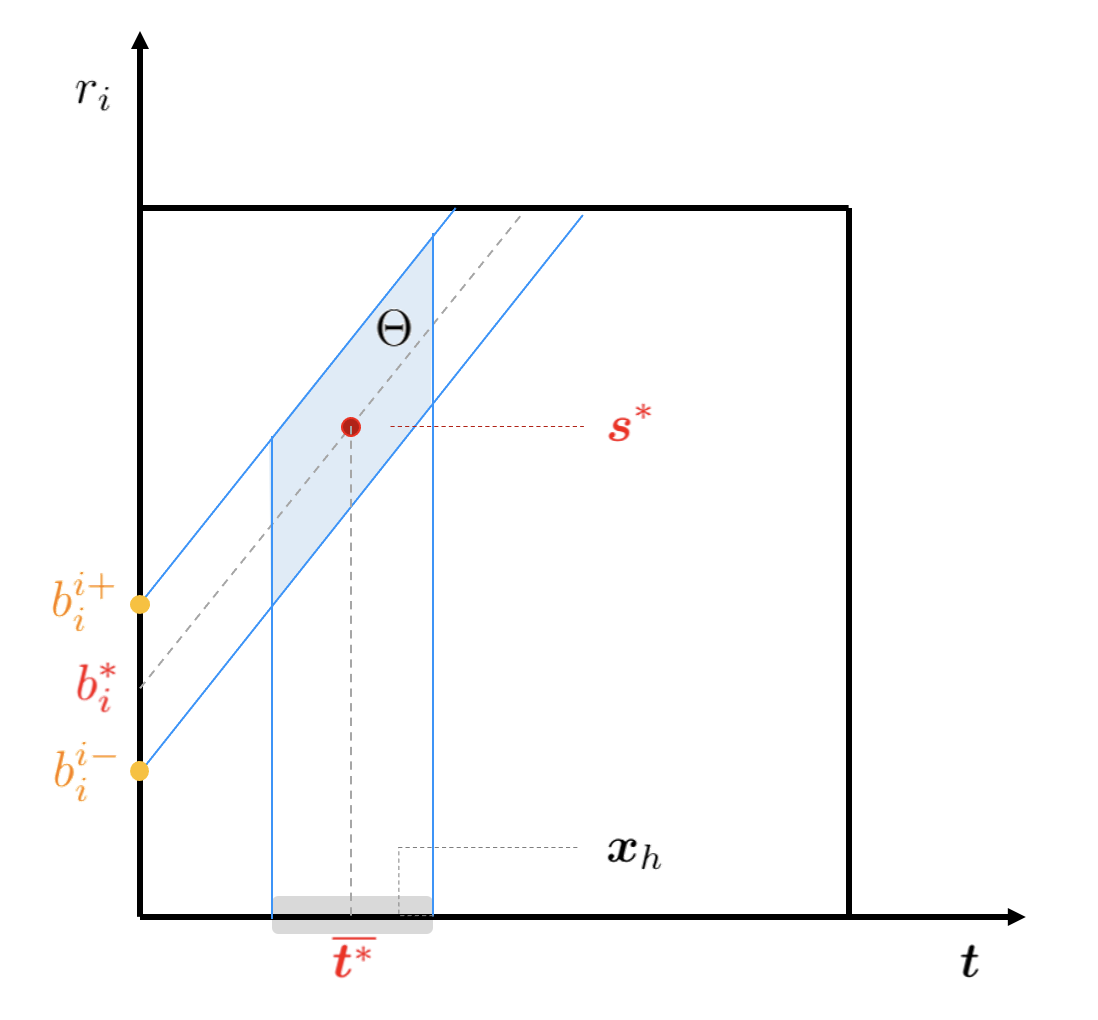}
	\caption{Process of determining $\mathcal{G}^i$. Here, the vertical axis exclusively considers the $i$-dimension in the $\bm{r}$-axis.}
    \label{fig:shrink}
\end{figure*}

Furthermore, according to $\overline{\bm{s}^*}=(\overline{\bm{t}^*},\overline{\bm{r}^*})^T:=\mathcal{M}(\bm{s}^*)$, the region $\mathcal{G}$ can be shrunk to $\bigcup\mathcal{G}^i$, where
\begin{equation*}
	\mathcal{G}^i:=\{\bm{x}\in \mathbb{R}^d:\bm{x}=\bm{x}^i+\overline{\bm{t}^*}+\bm{x}_h,\bm{x}^i\in\{\bm{x}^{i+},\bm{x}^{i-}\}, \bm{x}_h\in [-h/2,h/2)^d\},~~i=1,\cdots,n-d.
\end{equation*}
$\mathcal{C}(\bigcup\mathcal{G}^i)$ can directly enclose an interpolation element $\Theta$ in $\Omega$, and the size of interpolation element $\Theta$ fulfills
\begin{equation*}
   h-\varepsilon \leq\|\langle\Theta\rangle\|_{\infty}\leq h+\varepsilon, 
\end{equation*}
where $h$ is the target value of $\langle\Theta\rangle$ in each dimensional, and $\varepsilon$ determines the deviation range of $\langle\Theta\rangle$ in the actual selection. 
In concrete implementation, we usually choose $\varepsilon<h/10$ to ensure that $\langle\Theta\rangle$ is mainly controlled by $h$ and hardly affected by $\varepsilon$.
To facilitate a clearer comprehension of this step, \Cref{fig:shrink} presents the process of determining $\mathcal{G}^i$ in the $\bm{t}\times r_i$ plane. 
Here, the vertical axis exclusively considers the $i$-dimension in the $\bm{r}$-axis. 

\textbf{Interpolation recovery.}
To approximate $f(\bm{x}^*)$ by the $k$-degree Lagrange interpolation formula (LIF-$k$), we first select $K = (k+1)^n$ interpolation nodes $\{\bm{\tau}_i\in\mathbb{R}^d:i=1,\cdots,K\}$ in $\bigcup\mathcal{G}^i$. 
These nodes correspond to a point set $\{\mathcal{C}(\bm{\tau}_i)\}^K_{i=1}$, which forms an interpolation element $\Theta$ in $\Omega$. 
Before interpolation, affine transforming the irregular interpolation element into a rectangular shape is necessary.
Then the corresponding interpolation basis functions $\{\phi_i:\mathbb{R}^n\to\mathbb{R}:i=1,\cdots,K\}$ are determined.
Finally, combining \Cref{thm:relation},
\begin{equation*}
	\begin{aligned}
		f(\bm{x}^*)\approx\Pi_k F(\overline{\bm{s}^*})=\sum_{i=1}^{K}F(\mathcal{C}(\bm{\tau}_i))\phi_i(\overline{\bm{s}^*}) 
		=\sum_{i=1}^{K} f(\bm{\tau}_i)\phi_i(\overline{\bm{s}^*}),
	\end{aligned}
\end{equation*}
where $\Pi_k$ is the LIF-$k$ transformation.
According to \Cref{thm:dens_fill}, as the region $\mathcal{G}$ expands, $\mathcal{C}(\mathcal{G})$ gradually becomes dense in $\Omega$, and the approximation accuracy of $\Pi_k F(\overline{\bm{s}^*})$ is improved accordingly. 

\begin{remark}
    The LIF used in the FPR method can be replaced by other interpolation methods, and the shape of interpolation element can also be changed. One can make appropriate adjustments according to specific problems in the FPR framework. 
\end{remark}

\Cref{alg:FPR} summarizes the above steps of FPR method.

\begin{algorithm}[!pbht]
	\caption{FPR method}
	\label{alg:FPR}
	\begin{algorithmic}[1]
		\REQUIRE 
		projection matrix $\bm{P}\in \mathbb{P}^{d\times n}$,  target point $\bm{x}^*\in \mathbb{R}^d$
		\STATE Select region $\mathcal{G}\subset\mathbb{R}^d$
		\STATE Shrink region $\mathcal{G}$ to $\bigcup\mathcal{G}^i$
		\STATE Interpolation recovery
		\begin{itemize}
			\item Select interpolation nodes in  $\bigcup\mathcal{G}^i$
		    \item Use affine transformation to obtain rectangular interpolation element
		    \item Calculate interpolation basis functions
		    \item Calculate interpolation result  $\Pi_k F(\overline{\bm{s}^*})$
		\end{itemize}
	\end{algorithmic}
\end{algorithm}

\subsection{How to select region $\mathcal{G}$}\label{subsec:rapid_filling}

\begin{figure*}[!hbpt]
    \centering
	\subfigure[$\mathcal{C}_1(\mathcal{G})$ in $\mathbb{T}^n$]{\includegraphics[width=4.0cm]{./figures/torus_2.png}}
	\hspace{10mm}
	\subfigure[$\mathcal{C}_2(\mathcal{G})$ in $\mathbb{T}^n$]{\includegraphics[width=4.0cm]{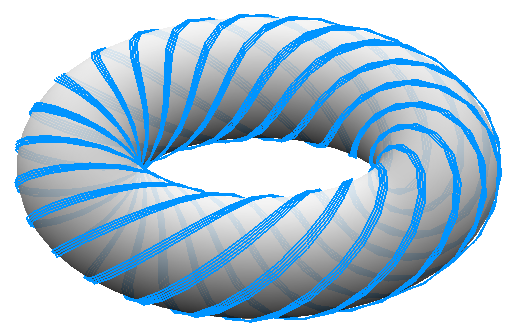}}	
	\caption{$\mathcal{C}_i(\mathcal{G})$ in $\mathbb{T}^n$, $i=1,2$, where $\mathcal{G}=[0, 30)$, $\mathcal{C}_1$ and $\mathcal{C}_2$ are determined by projection matrices $\bm{P}_1=(1,\sqrt{2})$ and $\bm{P}_2=(1,\pi)$, respectively.}
    \label{fig:torus_compare}
\end{figure*}

In the process of selecting region $\mathcal{G}$, 
we have observed that, in many cases, relatively small region $\mathcal{G}$ can result in a uniform distribution of $\mathcal{C}(\mathcal{G})$ in $\Omega$. 
Nevertheless, the distinct arithmetic properties of the irrational numbers in the projection matrices can influence the degree of uniformity in the distribution. 
For example, consider two different maps $\mathcal{C}_1$ and $\mathcal{C}_2$, which are determined by projection matrices $\bm{P}_1=(1,\sqrt{2})$ and $\bm{P}_2=(1,\pi)$, respectively. 
When $\mathcal{G}=[0, 30)$, the distribution of $\mathcal{C}_1(\mathcal{G})$ appears to be more uniform than that of $\mathcal{C}_2(\mathcal{G})$, as depicted in \Cref{fig:torus_compare}. 
Subsequently, we delve into the distribution characteristics of $\mathcal{C}(\mathcal{G})$ in two distinct categories: badly approximable systems and good approximable systems. Furthermore, we provide the respective selection strategies for the region $\mathcal{G}$ within each category. 

\subsubsection{Badly approximable systems}

As proved in \Cref{thm:dens_fill}, the distribution of $\mathcal{C}(\mathcal{G})$ in $\Omega$ is characterized by the distribution of $\{\bm{r}=\bm{A}\bm{t}:\bm{t}\in \mathcal{G}\cap \mathbb{Z}^d\}/\mathbb{Z}^{n-d}$ in $[0,1)^{n-d}$. 
When $A\in \mathbb{R}^{(n-d)\times d}$ is a badly approximable matrix, it has the rapid filling property mentioned in \Cref{rem:rapid_fill}. 
Specifically, given an $\varepsilon>0$, there is a least region $\mathcal{G}$ with
\begin{equation}
  (1/2)\varepsilon^{-(n-d)}\leq\langle\mathcal{G}(\varepsilon)\rangle_i\leq C_i\varepsilon^{-(n-d)},~~i=1,\cdots,d,  
\end{equation}
such that for each $\bm{r}_0\in [0,1)^{n-d}$, 
the $(n-d)$-dimensional Diophantine system
\begin{equation*}
        \|\bm{A}\bm{t}/\mathbb{Z}^{n-d}-\bm{r}_0\|_{\infty}\leq \varepsilon.
\end{equation*}
has a solution $\bm{t}\in\mathcal{G}(\varepsilon)\cap\mathbb{Z}^{d}$.
Thus, $\{\bm{r}=\bm{A}\bm{t}:\bm{t}\in\mathcal{G}(\varepsilon)\cap\mathbb{Z}^d\}/\mathbb{Z}^{n-d}$ is dense in $[0,1)^{n-d}$, and then $\mathcal{C}(\mathcal{G}(\varepsilon))$ is dense in $\Omega$. 

The rapid filling property offers an efficient way to select region $\mathcal{G}$ for badly approximable systems. 
Actually, $\varepsilon$ measures the \textit{filling precision}, and $\mathcal{G}(\varepsilon)$ is the \textit{filling cost} of achieving $\varepsilon$. 
For instance, we consider a badly approximable irrational number $\alpha=\sqrt{2}$.
\Cref{tab:region} shows the $\langle\mathcal{G}(\varepsilon)\rangle$ required for different filling precision $\varepsilon$. 
Specifically, given an $\varepsilon>0$, we traverse through $\mathbb{N}_+$ to search for a minimum point $t^*$ such that
\begin{equation*}
        \|\sqrt{2}t^*/\mathbb{Z}-0\|_{\infty}\leq \varepsilon
\end{equation*}
holds, and then $\langle\mathcal{G}(\varepsilon)\rangle=t^*$.
Results demonstrate an inverse relationship between $\varepsilon$ and $\langle\mathcal{G}(\varepsilon)\rangle$. 

\begin{table}[h]
    \centering
        \caption{Required $\langle\mathcal{G}(\varepsilon)\rangle$ for the badly approximable irrational number $\alpha=\sqrt{2}$ under different filling precision $\varepsilon$.}
    \label{tab:region}
    \begin{tabular}{|c|c|c|c|c|c|}
    \hline
	$\varepsilon$ & 1.0e-02 &  1.0e-03 & 1.0e-04 & 1.0e-05 & 1.0e-06\\
	\hline
	$\langle\mathcal{G}(\varepsilon)\rangle$ & 9.9e+01 & 9.85e+02 & 5.74e+03 & 1.14e+05 & 1.13e+06\\
    \hline
    \end{tabular}
\end{table}

Therefore, when using FPR method to recover the quasiperiodic functions related to badly approximable systems, we can directly select the region $\mathcal{G}$ to reach the upper bound of $\langle\mathcal{G}(\varepsilon)\rangle$, i.e., $\mathcal{G}$ fulfills
\begin{equation}\label{eq:bad_region}
    \langle\mathcal{G}\rangle_i= C_i\varepsilon^{-(n-d)},~~i=1,\cdots,d,
\end{equation}
where $C_i$ is the smallest integer upper bound of $\langle\mathcal{G}(\varepsilon)\rangle_i\varepsilon^{(n-d)}$.

\subsubsection{Good approximable systems}\label{subsub:good_appro}

For each good approximable system, there is no uniform upper bound of $\langle\mathcal{G}(\varepsilon)\rangle$ that grows linearly with respect to $\varepsilon^{-(n-d)}$, indicating the absence of rapid filling property.    
Consequently, the determination of $\mathcal{G}$ under this category cannot proceed in the same manner as it does for badly approximable systems. 
For example, consider the good approximable system
\begin{equation}\label{eq:good_pi}
        \|\pi t/\mathbb{Z}-r_0\|_{\infty}\leq \varepsilon,~~\forall r_0\in \mathbb{T}, 
\end{equation}
where $\pi$ is a well-known transcendental number. 
\Cref{tab:approx_region} presents the $\langle\mathcal{G}(\varepsilon)\rangle$ of the good approximable system \cref{eq:good_pi} required for each given filling precision $\varepsilon$. $\langle\mathcal{G}(\varepsilon)\rangle$ shows an exponential-like growth behavior, without an inverse relationship between $\langle\mathcal{G}(\varepsilon)\rangle$ and $\varepsilon$. 
\begin{table}[!hbpt]
    \centering
        \caption{Required $\langle\mathcal{G}(\varepsilon)\rangle$ for the good approximable irrational number $\pi$ under different filling precision $\varepsilon$.}
    \label{tab:approx_region}
    \begin{tabular}{|c|c|c|c|c|c|c|c|c|c|c|}
    \hline
	$\varepsilon$ & 1.0e-03 & 1.0e-04 & 1.0e-05 & 1.0e-06 & 1.0e-07\\
	\hline
	$\langle\mathcal{G}(\varepsilon)\rangle$ & 2.05e+05 & 4.16e+05 & 6.25e+06 & 1.27e+08 & $>$1.0e+10\\
    \hline
    \end{tabular}
\end{table}

However, we can also propose a convenient and efficient approach of selecting $\mathcal{G}$ for the quasiperiodic function related to good approximable system.
According to \Cref{lem:lower_bound}, the lower bound $\langle\mathcal{G}(\varepsilon)\rangle_i\geq (1/2)\varepsilon^{-(n-d)}$, $i=1,\cdots,d$, holds for each Diophantine system. 
It allows us, for the good approximable system, to adopt the lower bound $\langle\mathcal{G}\rangle_i=(1/2)\varepsilon^{-(n-d)}$, $i=1,\cdots,d$ to determine $\mathcal{G}$, named the computable region. 
In this way, for a given $\varepsilon$ the size of $\mathcal{G}$ can be immediately determined without traversing, resulting in a significantly reduce of computational costs.

\subsection{Computational complexity analysis}

Here, we analyze the computational complexity of FPR method for recovering a required interpolation point $\bm{x}^*$.  
Assume that the size of interpolation element $\Theta\subset\Omega$ fulfils
\begin{equation*}
   h-\varepsilon \leq\|\langle\Theta\rangle\|_{\infty}\leq h+\varepsilon,
\end{equation*}
where $\varepsilon$ is the filling precision satisfying $\varepsilon<h/10$. Here we fix $\varepsilon=h/\tilde{c}$, where constant $\tilde{c}>10$. $N=1/h$ denotes the degree of freedom of spatial discrete nodes. Then
$\|\langle\mathcal{G}\rangle\|_{\infty}=C\varepsilon^{-(n-d)}=C(\tilde{c}N)^{n-d}$, where $C=\max_{i=1}^dC_i$. 
The traverse process in $\mathcal{G}\cap\mathbb{Z}^d$ can be performed simultaneously in each dimension,  
and computing $\|\bm{b}_i-\bm{b}^*\|_{\infty}$ consumes $O(n-d)$ times of subtraction. 
Hence, the computational complexity of obtaining the interpolation element of $\bm{x}^*$ is 
\begin{equation*}
    O\left(d(n-d)\|\langle\mathcal{G}\rangle\|_{\infty}\right)=O\left(d(n-d)C\varepsilon^{-(n-d)}\right)=O\left(N^{n-d}\right).
\end{equation*}

To use an $n$-dimensional LIF-$k$, $K=(k+1)^n$ interpolation nodes are required. 
Thus, the computational complexity of interpolation is $O(K)$. 
Since $K\ll N$, the computational complexity of the FPR method is at the level of $O\left(N^{n-d}\right)$. 
It is worth emphasizing that when using the FPR method for practical problems, the interpolation nodes can be predetermined, which means that the computational cost of using the FPR method could be almost negligible.

\subsection{Convergence analysis}
In this subsection, we provide the convergence analysis of the FPR method for recovering a target point within an interpolation element $\Theta\subset\Omega\subset\mathbb{R}^n$. 
To facilitate this analysis, we introduce the following multi-index notations.  
For a vector $\bm{\alpha}=(\alpha_i)^n_{i=1}\in \mathbb{R}^n$, 
$|\bm{\alpha}|:=\alpha_1+\cdots+\alpha_n$. 
$\mathcal{C}^{\bm{\alpha}}(\mathbb{R}^n)$ denotes the set of functions on $\mathbb{R}^n$ with $\alpha_i$-order continuous derivatives along the $i$-th coordinate direction. 
For a function $F(\bm{s})\in \mathcal{C}^{\bm{\alpha}}(\mathbb{R}^n)$, its $\bm{\alpha}$-order derivative is
$$D^{\bm{\alpha}}F(\bm{s}):= \dfrac{\partial^{|\bm{\alpha}|}F(\bm{s})}{\partial s_1^{\alpha_1}\cdots\partial s_n^{\alpha_n}}.$$
Define space
\begin{equation*}
    \mathcal{L}^2(\Omega):=\{F:\|F\|_{\mathcal{L}^2(\Omega)}< \infty\},
\end{equation*}
where
\begin{equation*}
    \|F\|_{\mathcal{L}^2(\Omega)}=\Bigg(\int_{\Omega}|F(\bm{s})|^2d\bm{s}\Bigg)^{1/2}.
\end{equation*}
For any integer $m\geq 0$, the Hilbert space on $\Omega$ is 
\begin{equation*}
    \mathcal{H}^m(\Omega):=\{F\in \mathcal{L}^2(\Omega):\|F\|_{m,\Omega}< \infty\},
\end{equation*}
where
\begin{equation*}
    \|F\|_{m,\Omega}=\left(\sum_{|\bm{\alpha}|\leq m}\|D^{\bm{\alpha}}F\|^2_{\mathcal{L}^2(\Omega)}\right)^{1/2}.
\end{equation*}
The semi-norm of $\mathcal{H}^m(\Omega)$ is defined as 
\begin{equation*}
    |F|_{m,\Omega}=\left(\sum_{|\bm{\alpha}|= m}\|D^{\bm{\alpha}}F\|^2_{\mathcal{L}^2(\Omega)}\right)^{1/2}.
\end{equation*}

\begin{definition}\label{def:affine}
    Two elements $\Theta$ and $\hat{\Theta}$ are affine equivalent, if there is an invertible affine transformation
    \begin{equation*}
        \begin{aligned}
        &\mathcal{A}: &\hat{\Theta}&\to \Theta,\\
        &&\hat{\bm{s}}&\mapsto \bm{s}=\bm{B}\hat{\bm{s}}+\bm{b},
        \end{aligned}
    \end{equation*}
    where $\bm{B}\in\mathbb{R}^{n\times n}$ is an invertible matrix and $\bm{b}\in \mathbb{R}^n$.
\end{definition}

\begin{lemma}\cite{scott1990finite}
    If two elements $\Theta$ and $\hat{\Theta}$ are affine equivalent, for any $F\in \mathcal{H}^m(\Theta)$, let 
    \begin{equation*}
        \hat{F}:F\circ \mathcal{A},
    \end{equation*}
    then $\hat{F}\in\mathcal{H}^m(\hat{\Theta})$ and there is a positive constant $C=C(m,n)$ such that
    \begin{equation*}
        \begin{aligned}
        |\hat{F}|_{m,\hat{\Theta}} &\leq C\|\bm{B}\|^m|\det(\bm{B})|^{-\frac{1}{2}}|F|_{m,\Theta},\\
        |F|_{m,\Theta}  &\leq C\|\bm{B}^{-1}\|^m|\det(\bm{B})|^{\frac{1}{2}}|\hat{F}|_{m,\hat{\Theta}}.
        \end{aligned}
    \end{equation*}
    Here, $\|\cdot\|$ is the Euclidean norm.
\end{lemma}

\begin{lemma}\cite{scott1990finite}
    The matrix $\bm{B}$ defined in \Cref{def:affine} satisfies
    \begin{equation*}
        \|\bm{B}\|\leq \frac{h_{\Theta}}{\hat{\rho}_{\hat{\Theta}}},~~\|\bm{B}^{-1}\|\leq \frac{\hat{h}_{\hat{\Theta}}}{\rho_{\Theta}},
    \end{equation*}
    where
    \begin{equation*}
        \begin{cases}
            h_{\Theta}=\rm{diam}\Theta,&\rho_{\Theta}=\sup\{\rm{diam}\Gamma:closed~ball~\Gamma\subset\Theta\},\\
            \hat{h}_{\hat{\Theta}}=\rm{diam}\hat{\Theta},&\hat{\rho}_{\hat{\Theta}}=\sup\{\rm{diam}\Gamma:closed~ball~ \Gamma\subset\hat{\Theta}\}.
        \end{cases}
    \end{equation*}
\end{lemma}

Note that the interpolation element $\Theta$ is a regular shape, i.e., there exists a constant $\kappa>0$ such that $h_{\Theta}/\rho_{\Theta}\leq \kappa$.
Denote $P_k(\hat{\Theta})$ as the space consisting of polynomials of degree $k$ or less in $\hat{\Theta}$, 
and $\hat{\Pi}_k$ as the LIF-$k$ operator over $\hat{\Theta}$. We have the convergence result of the FPR method.

\begin{theorem}
    If $0\leq m\leq k+1$, then for $\hat{\Pi}_k\in\mathcal{L}(\mathcal{H}^{k+1}(\hat{\Theta});\mathcal{H}^m(\hat{\Theta}))$, $\hat{\Pi}_k\hat{F}=\hat{F}$ for each $\hat{F}\in P_k(\hat{\Theta})$.
    Here, $\mathcal{L}(A; B)$ denotes the set of linear operators from space $A$ to space $B$.
    Over the interpolation element $\Theta~(\text{affine equivalent with}~\hat{\Theta})$, operator $\Pi_k$ fulfils
    \begin{equation*}
        (\Pi_kF)^{\hat{}}=\hat{\Pi}_k\hat{F},~~ \hat{F}\in P_k(\hat{\Theta}).
    \end{equation*}
    Then there is a constant $C$ such that
    \begin{equation*}
        \|F-\Pi_kF\|_{m,\Theta}\leq Ch_{\Theta}^{k+1-m}|F|_{k+1,\Theta},~~F\in \mathcal{H}^{k+1}(\Theta).
    \end{equation*}
    \label{thm:conv}
\end{theorem}
\begin{proof}
    Since there is no error in the process of selecting finite points, the only error in the FPR method comes from the interpolation error. According to the interpolation error analysis in \cite{scott1990finite}, the above conclusion can be established.
\end{proof}

\section{Numerical experiments}\label{sec:num}
In this section, we present two classes of quasiperiodic systems to show the performance of FPR method.
One class is analytical quasiperiodic functions, another is a piecewise constant quasicrystal.
    
We use the $\ell^{\infty}$-norm to measure the error between numerical result and exact result, denoted as $e(\langle\Theta\rangle)$ with respect to the interpolation element $\Theta$. 
Then we can estimate the order of accuracy by calculating the logarithmic ratio of errors between two successive refinements
\begin{equation*}
    \text{Order}=\log_{2}\left(\frac{e(2\langle\Theta\rangle)}{e(\langle\Theta\rangle)}\right).
\end{equation*}

To apply the FPR method for recovering $d$-dimensional quasiperiodic functions related to badly approximable systems, given a filling precision $\varepsilon$, 
we select the region $\mathcal{G}$ given by
\begin{equation*}
    \mathcal{G}=[0,C_1\varepsilon^{-(n-d)})\times \cdots \times[0,C_d\varepsilon^{-(n-d)}),
\end{equation*}
where $n$ is the dimension of superspace and $C_i~(i=1,\cdots,d)$ is determined by
\begin{equation*}
    C_i=\max_{j=0,\cdots,4}\langle\mathcal{G}(\varepsilon_j)\rangle_i\varepsilon^{n-d}_j,~~\varepsilon_j=\frac{1}{10^{2+j}}.
\end{equation*}
Since $j$ is small, determining constant $C_i~(i=1,\cdots,d)$ in this way consumes almost no computational cost.

As discussed in \Cref{subsub:good_appro}, for the quasiperiodic functions related to good approximable systems, 
we select the computable region
\begin{equation*}
    \mathcal{G}=[0,0.5\varepsilon^{-(n-d)})^d.
\end{equation*}

\subsection{Smooth quasiperiodic functions}
In this subsection, we present four different smooth quasiperiodic functions to demonstrate the accuracy and efficiency of FPR method. The first three examples are related to badly approximable systems, and the last one is related to a good approximable system.
\begin{example}\label{exmp:1Dbad}
	Consider a 1D quasiperiodic function
    \begin{equation}\label{eq:2to1}
        f(x)=\cos x+\cos\sqrt{2}x,\quad x\in \mathbb{R}.
    \end{equation}
\end{example} 

$f(x)$ can be embedded into the 2D periodic function $F(\bm{s})=\cos\bm{s}$,  $\bm{s}=(s_1,s_2)^T\in \mathbb{R}^2$ using the projection matrix $\bm{P}=(1,\sqrt{2})$. 
The unit cell of $F(\bm{s})$ is $\Omega=[0,2\pi)^2$. 
Since $\sqrt{2}$ is a badly approximable irrational number, the corresponding Diophantine system has the rapid filling property. 
\Cref{tab:2to1_region} shows the inverse relationship between the filling precision $\varepsilon$ and the size of region $\mathcal{G}(\varepsilon)$. 
Therefore, selecting a relatively smaller region $\mathcal{G}$ for FPR method is sufficient to recover the global information of quasiperiodic function \cref{eq:2to1}.

\begin{table}[!hbpt]
    \centering
        \caption{Required $\langle\mathcal{G}(\varepsilon)\rangle$ with different filling precision $\varepsilon$ for projection matrix $\bm{P}=(1,\sqrt{2})$.}
    \label{tab:2to1_region}
    \begin{tabular}{|c|c|c|c|c|c|}
    \hline
	$\varepsilon$ & 1.0e-02 &  1.0e-03 & 1.0e-04 & 1.0e-05 & 1.0e-06\\
	\hline
	$\langle\mathcal{G}(\varepsilon)\rangle/2\pi$ & 9.90e+01 & 9.85e+02 & 5.74e+03 & 1.14e+05 & 1.13e+06\\
    \hline
    \end{tabular}
\end{table}

\begin{figure*}[!hbpt]
    \centering
	\subfigure[$\langle\Theta\rangle=(0.4,0.3)$]{\label{fig:linear_a}\includegraphics[width=5cm]{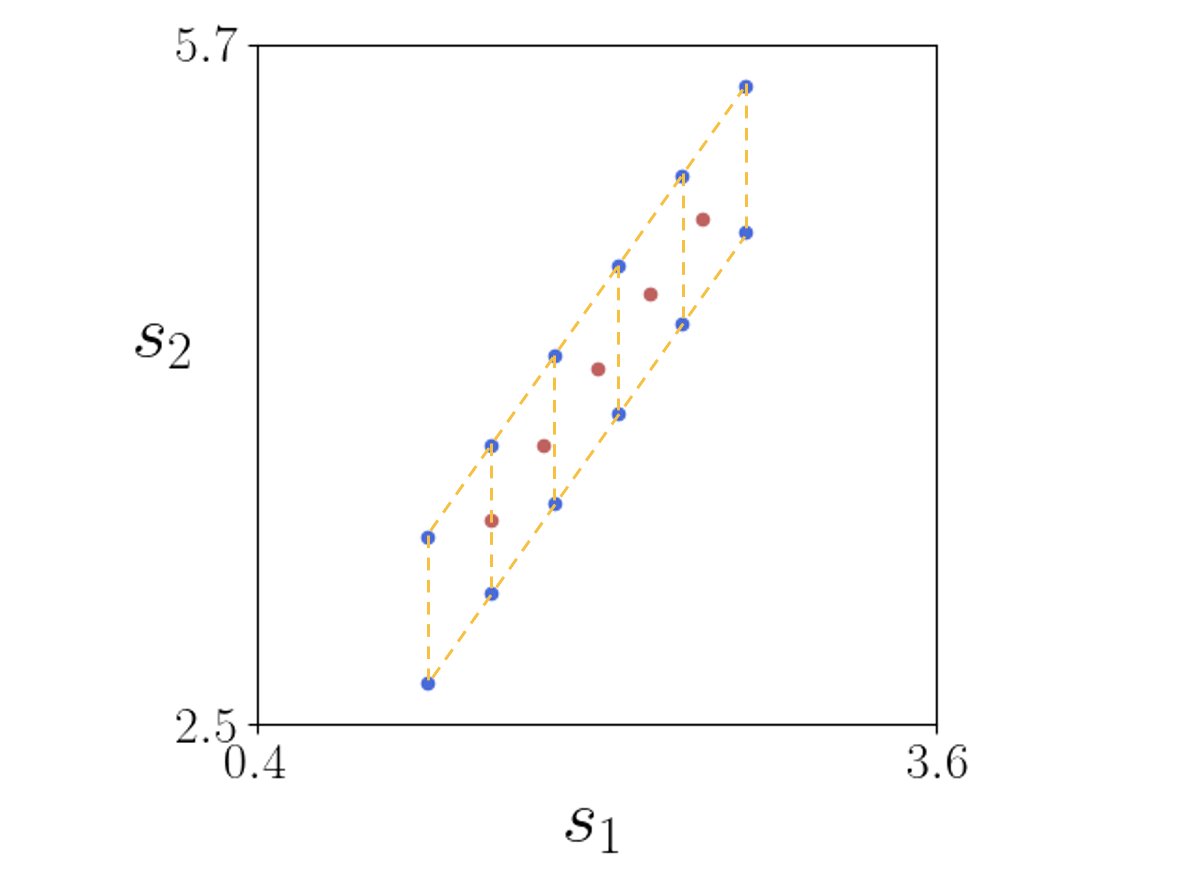}}
	\hspace{5mm}
  	\subfigure[$\langle\Theta\rangle=(0.2,0.15)$]{\label{fig:linear_b}\includegraphics[width=5cm]{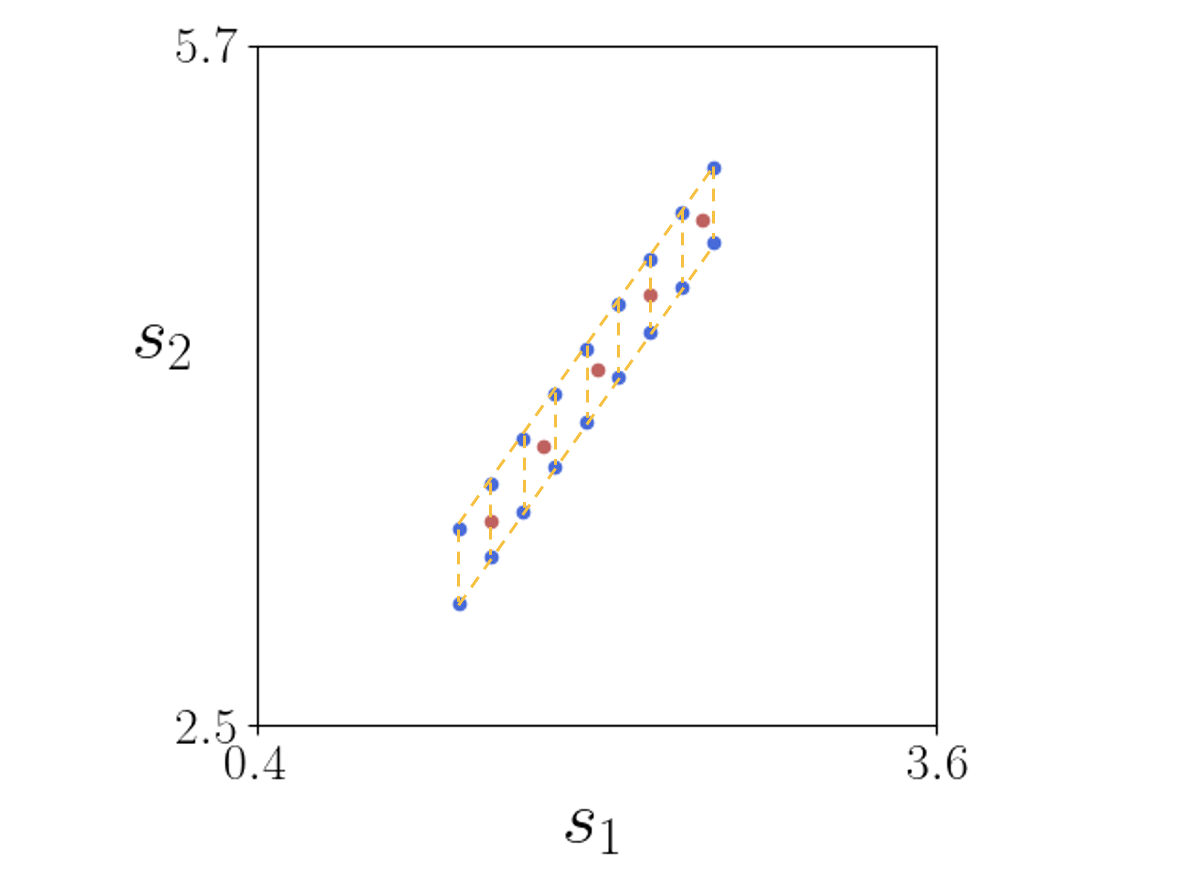}}
	\caption{Node selection for FPR method with 2D LIF-1 when solving the quasiperiodic function \cref{eq:2to1}, $x\in [6284, 6286)$, with different sizes of interpolation elements $\Theta$.}
    \label{fig:linear}
\end{figure*}

First, we use the FPR method with 2D LIF-1 in $\mathcal{G}=[0,8\varepsilon^{-1})$ to recover \cref{eq:2to1} when $x\in[6284,6286)$. 
\Cref{fig:linear} shows the node selection of FPR method with two sizes of interpolation elements. 
For each target point $x^*\in [6284, 6286)$, its image $\overline{\bm{s}^*}$ under the combination map $\mathcal{C}$ is marked with red dot. 
We have circled each interpolation element $\Theta$ that contains the target point using yellow dashed lines. 
Four blue dots in each $\Theta$ are used as interpolation nodes. 
In \Cref{fig:linear_a}, the size of interpolation elements is $\langle\Theta\rangle=(0.4,0.3)$.
After one refinement, the size of interpolation element becomes $\langle\Theta\rangle=(0.2,0.15)$, as shown in \Cref{fig:linear_b}.

\begin{table}[!hbpt]
    \centering
    \caption{Error of FPR method when recovering quasiperiodic function \cref{eq:2to1} for different interpolation element sizes $\langle\Theta\rangle$.}
    \label{tab:example1_error}
    \resizebox{\textwidth}{!}{
    \begin{tabular}{|c|c|c|c|c|c|c|c|c|}
    \hline
    \multicolumn{3}{|c|}{FPR with 2D LIF-1} & \multicolumn{3}{c|}{FPR with 2D LIF-3} & \multicolumn{3}{c|}{FPR with 2D LIF-5}\\
    \hline
    $\langle\Theta\rangle$ & Error & Order & $\langle\Theta\rangle$ & Error & Order & $\langle\Theta\rangle$ & Error & Order \\
	\hline
	(0.4, 0.3) & 2.8035e-02& & (0.8, 0.3) & 2.6990e-03 & & (1.2, 0.3) & 9.6976e-05 &\\
	(0.2, 0.15) & 7.0599e-03 & 1.99 & (0.4, 0.15) & 1.7661e-04 & 3.93 & (0.6, 0.15) & 1.6192e-06 & 5.90\\
	(0.1, 0.075) & 1.7681e-03 & 2.00 & (0.2, 0.075) & 1.1204e-05 & 3.98 & (0.3, 0.075) & 2.5909e-08 & 5.97\\
	(0.05, 0.0375) & 4.4219e-04 & 2.00 & (0.1, 0.0375) & 7.0414e-07 & 3.99 & (0.15, 0.0375) & 4.0196e-10 & 6.01\\
    \hline
    \end{tabular}
    }
\end{table} 

\Cref{tab:example1_error} records the errors and accuracy orders obtained by using the FPR method with 2D LIF-1, LIF-3, and LIF-5 to recover the quasiperiodic function \cref{eq:2to1} when $x\in[6284,6286)$. 
As \Cref{thm:conv} predicts, using higher-order interpolation methods can result in better approximation accuracy.

\begin{figure*}[!hbpt]
    \centering
	\includegraphics[width=0.8\textwidth]{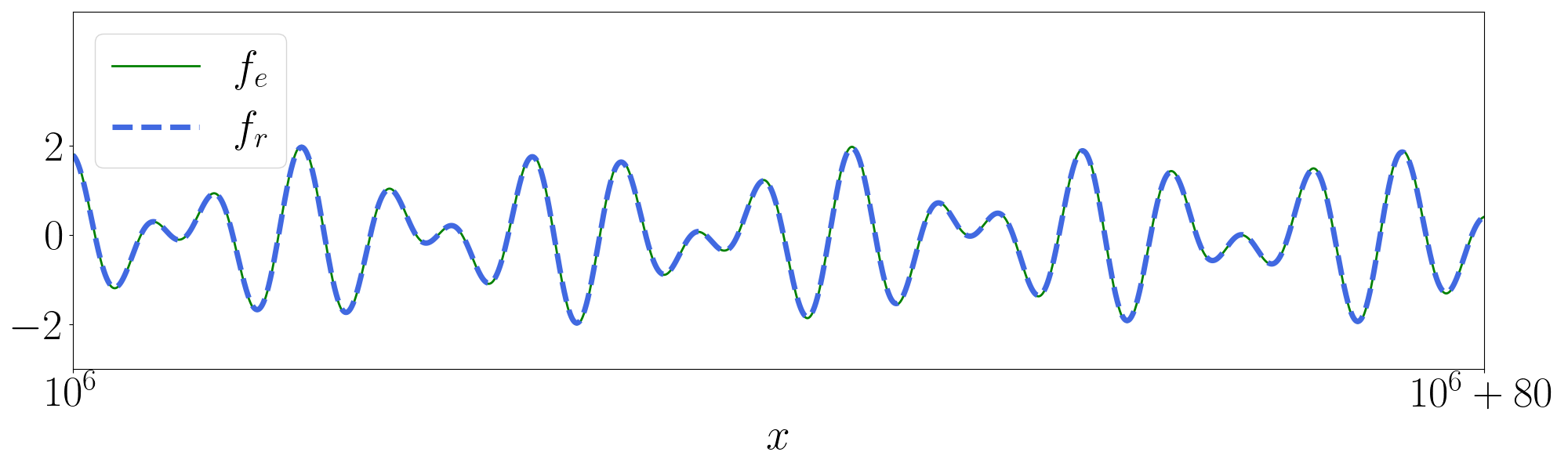}
	\caption{Comparison of the exact value $f_e$ and the recovered result $f_r$ by FPR method with 2D LIF-1 for quasiperiodic function \cref{eq:2to1} when $x\in[10^6,10^6+80)$.}
    \label{fig:2to1_global}
\end{figure*}

Next, we further present the power of FPR method on recovering the global quasiperiodic function with finite points. 
Here, we select the region $\mathcal{G}=[0,8000)$ and the size of the interpolation elements $\langle\Theta\rangle=(0.4,0.3)$. 
In this case, we only need 320 interpolation nodes to recover the global information of $f(x)$.
\Cref{fig:2to1_global} compares the recovered result by FPR method and the exact value of \cref{eq:2to1} when $x\in[10^6, 10^6+80)$. The error between the two is 3.1147e-02.
 
\begin{example}
    Consider a 1D quasiperiodic function
    \begin{equation}\label{eq:3to1}
        f(x)=\cos x+\cos\sqrt{2}x+\cos\sqrt{3}x,\quad x\in \mathbb{R}.
    \end{equation}
\end{example}

Compared with \Cref{exmp:1Dbad} only containing one irrational frequency $\sqrt{2}$, this quasiperiodic function $f(x)$ has two irrational frequencies $\sqrt{2}$ and $\sqrt{3}$. Correspondingly, $f(x)$ should embed into the 3D periodic function $F(\bm{s})=\cos\bm{s}$ with unit cell $\Omega=[0,2\pi)^3$, through the projection matrix $\bm{P}=(1,\sqrt{2},\sqrt{3})$. 
Due to $\sqrt{2}$ and $\sqrt{3}$ are both  quadratic irrational numbers, $f(x)$ is related to a badly approximable system. 
There is an inverse relationship between the size of $\mathcal{G}(\varepsilon)$ and the filling precision $\varepsilon$ as shown in \Cref{tab:3to1_region}. 
Concretely, the size of $\mathcal{G}(\varepsilon)$ is the maximum size of $\mathcal{G}(\varepsilon)$ produced by   $\bm{P}_1=(1,\sqrt{2})$ and $\bm{P}_2=(1,\sqrt{3})$ for a given filling precision $\varepsilon$.
From this example, one can find that, as the dimension of superspace increases, however, the required interpolation points in $\mathcal{G}$ still belong to $\mathbb{R}$. 
It demonstrates that the FPR method is a no-lift algorithm. 

\begin{table}[!hbpt]
    \centering      
    \caption{Filling precision $\varepsilon$ and required $\langle\mathcal{G}(\varepsilon)\rangle$ for different projection matrices.}
    \label{tab:3to1_region}
    \resizebox{\textwidth}{!}{
    \begin{tabular}{|c|c|c|c|c|c|}
    \hline
	$\varepsilon$ & 1.0e-02 &  1.0e-03 & 1.0e-04 & 1.0e-05 & 1.0e-06\\
	\hline
	$\langle\mathcal{G}(\varepsilon)\rangle/2\pi$ for $\bm{P}_1=(1,\sqrt{2})$& 9.90e+01 & 9.85e+02 & 5.74e+03 & 1.14e+05 & 1.13e+06\\
	\hline
	$\langle\mathcal{G}(\varepsilon)\rangle/2\pi$ for $\bm{P}_2=(1,\sqrt{3})$ & 9.70e+01 & 1.35e+03 & 7.95e+03 & 1.10e+05 & 7.99e+05 \\
	\hline
	$\langle\mathcal{G}(\varepsilon)\rangle/2\pi$ for $\bm{P}=(1,\sqrt{2},\sqrt{3})$ & 9.90e+01 & 1.35e+03 & 7.95e+03 & 1.14e+05 & 1.13e+06\\
    \hline
    \end{tabular}}
\end{table}

\begin{figure*}[!hbpt]
    \centering
	\subfigure[$\langle\Theta\rangle=(0.4,0.4,0.3)$]{\label{fig:3to1_a}\includegraphics[width=5cm]{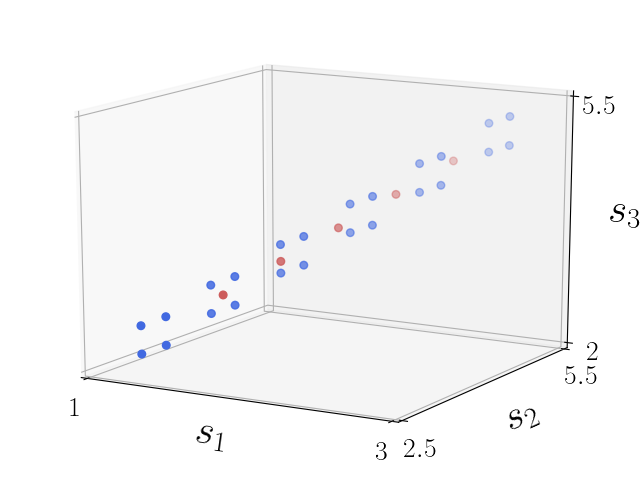}}
	\hspace{5mm}
	\subfigure[$\langle\Theta\rangle=(0.2,0.2,0.15)$]{\label{fig:3to1_b}\includegraphics[width=5cm]{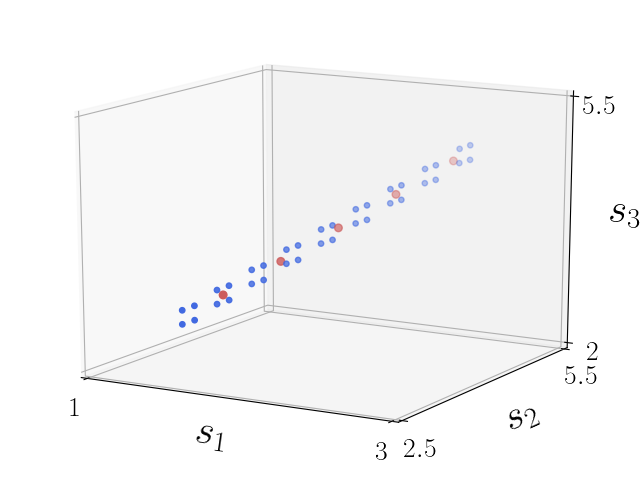}}
	\caption{Node selection for 3D LIF-1 FPR method with different interpolation element $\Theta$ to recover quasiperiodic function \cref{eq:3to1} on the interval $[6284, 6286)$.}
    \label{fig:3to1}
\end{figure*}

We then employ the FPR method with 3D LIF-1 in the region $\mathcal{G}=[0,9\varepsilon^{-1})$, to recover \cref{eq:3to1} when $x\in [6284, 6286)$. 
The 3D LIF-1 requires eight interpolation nodes within each interpolation element $\Theta$. 
\Cref{fig:3to1} shows the node selection for FPR method with two sizes of interpolation elements. 
\Cref{tab:example2_error} lists the error and accuracy order, consistent with the prediction by \Cref{thm:conv}.

\begin{table}[!hbpt]
    \centering
        \caption{Error of 3D LIF-1 FPR method for recovering quasiperiodic function \cref{eq:3to1} with interpolation element sizes $\langle\Theta\rangle$.}
    \label{tab:example2_error}
    \begin{tabular}{|c|c|c|}
    \hline
    $\langle\Theta\rangle$ & Error & Order \\
	\hline
	(0.4, 0.4, 0.3) & 1.1214e-01 & \\
	(0.2, 0.2, 0.15) & 2.7964e-02 & 2.00 \\
	(0.1, 0.1, 0.075) & 6.9445e-03 & 2.01 \\
	(0.05, 0.05, 0.0375) & 1.7093e-03 & 2.02 \\
    \hline
    \end{tabular}
\end{table}

\begin{example}
    Consider a 2D quasiperiodic function
    \begin{equation}\label{eq:3to2}
    f(x, y)=\cos x+\cos\sqrt{2}x+\cos y,\quad (x,y)\in \mathbb{R}^2.
    \end{equation}
\end{example}

There exists a projection matrix 
    \begin{equation}\label{eq:3to2_pmatrix}
        \bm{P}=\begin{bmatrix}
        1 & \sqrt{2} & 0 \\
        0 & 0 & 1\\
        \end{bmatrix}
    \end{equation}
such that $f(x,y)=F(\bm{P}^T(x,y)^T)$, 
where $F(\bm{s})=\cos\bm{s}$, $\bm{s}\in \mathbb{R}^3$, is the parent function of $f(x, y)$. 
The unit cell is $\Omega=[0,2\pi)^3$.
\Cref{tab:3to2_region} shows the inverse relationship between the filling precision $\varepsilon$ and $\langle\mathcal{G}\rangle_1$. 
Since $f(x,y)$ is periodic in the $y$ direction, it is sufficient to take $\langle\mathcal{G}\rangle_2$ as $2\pi$. 
Therefore, given the filling precision $\varepsilon$, we select $\mathcal{G}=[0, 8\varepsilon^{-1})\times[0,2\pi)$.

\begin{table}[!hbpt]
    \centering
        \caption{Required $\langle\mathcal{G}(\varepsilon)\rangle$ of different filling precision $\varepsilon$ for projection matrix \cref{eq:3to2_pmatrix}.}
    \label{tab:3to2_region}
    \begin{tabular}{|c|c|c|c|c|c|}
    \hline
	$\varepsilon$ & 1.0e-02 &  1.0e-03 & 1.0e-04 & 1.0e-05 & 1.0e-06\\
	\hline
	$\langle\mathcal{G}(\varepsilon)\rangle_1/2\pi$ & 9.90e+01 & 9.85e+02 & 5.74e+03 & 1.14e+05 & 1.13e+06\\
	\hline
	$\langle\mathcal{G}(\varepsilon)\rangle_2/2\pi$ & 1.0 & 1.0 & 1.0 & 1.0 & 1.0\\
    \hline
    \end{tabular}
\end{table}

\begin{figure*}[h]
    \centering
	\subfigure[$\langle\Theta\rangle=(0.8,0.3,0.3)$]{\label{fig:3to2_a}\includegraphics[width=5cm]{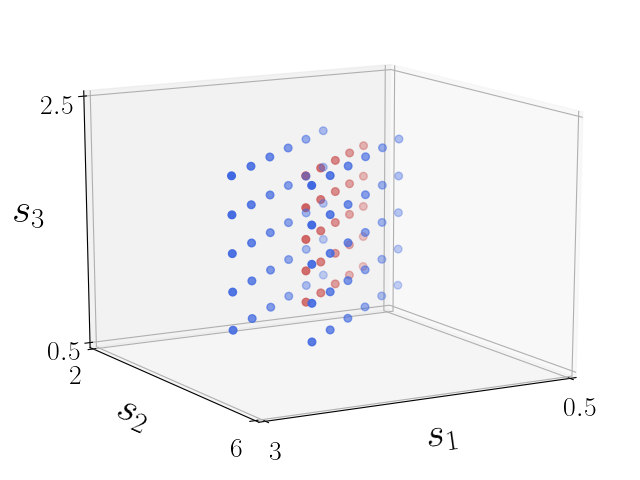}}
	\hspace{5mm}
	\subfigure[$\langle\Theta\rangle=(0.4,0.15,0.15)$]{\label{fig:3to2_b}\includegraphics[width=5cm]{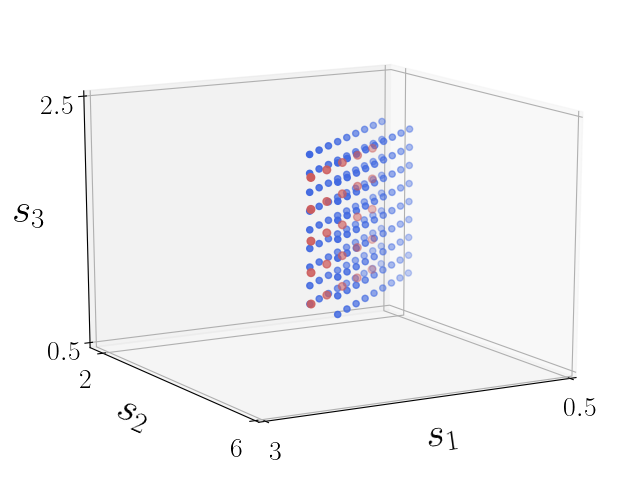}}
	\caption{Node selection of FPR method with 3D LIF-1 when recovering quasiperiodic function \cref{eq:3to2}, $(x, y)\in [6284, 6286)^2$, with different sizes of interpolation elements $\Theta$.}
    \label{fig:3to2}
\end{figure*}

We use FPR method with 3D LIF-1 to recover \cref{eq:3to2} when $(x,y)^T\in [6284, 6286)^2$. 
The 3D LIF-1 requires eight interpolation nodes in each interpolation element $\Theta$.  
\Cref{fig:3to2} shows the node selection for FPR method with different interpolation elements.\Cref{tab:example3_error} records the error and accuracy order, which is consistent with theoretical results. 

\begin{table}[!hbpt]
    \centering
        \caption{Error for FPR method with 3D LIF-1 when solving quasiperiodic function \cref{eq:3to2} for interpolation element sizes $\langle\Theta\rangle$.}
    \label{tab:example3_error}
    \begin{tabular}{|c|c|c|}
    \hline
    $\langle\Theta\rangle$ & Error & Order \\
	\hline
	(0.8, 0.3, 0.3) & 1.1654e-01 & \\
	(0.4, 0.15, 0.15) & 2.6666e-02 & 2.13 \\
	(0.2, 0.075, 0.075) & 5.8759e-03 & 2.18 \\
    \hline
    \end{tabular}
\end{table}

\begin{example}\label{exa:tran}
    Consider a 1D quasiperiodic function with a transcendental frequency $\pi$
    \begin{equation}\label{eq:approx}
        f(x)=\cos x+\cos\pi x,\quad x\in \mathbb{R}.
    \end{equation}
\end{example}

\begin{table}[h]
    \centering
        \caption{Required sizes $\langle\mathcal{G}\rangle$ of computable region $\mathcal{G}$ for different filling precision $\varepsilon$.}
    \label{tab:suboptimal_region}
    \begin{tabular}{|c|c|c|c|c|c|c|c|c|c|c|}
    \hline
	$\varepsilon$ & 1.0e-02 & 1.0e-03 & 1.0e-04 & 1.0e-05 & 1.0e-06\\
	\hline
	$\langle\mathcal{G}\rangle/2\pi$ & 5.0e+01 & 5.0e+02 & 5.0e+03 & 5.0e+04 & 5.0e+05 \\
    \hline
    \end{tabular}
\end{table} 

The projection matrix is $\bm{P}=(1,\pi)$. 
Since $\pi$ is a good approximable number, 
we can directly select finite points in the computable region $\mathcal{G}=[0,0.5\varepsilon^{-1})$ as discussed in \Cref{subsub:good_appro}.
The sizes of required computable region $\mathcal{G}$ with different filling precision $\varepsilon$ are shown in \Cref{tab:suboptimal_region}.

We apply the FPR method with 2D LIF-1 to recover \cref{eq:approx} when $x\in[6284,6286)$. 
As demonstrated in \Cref{exmp:1Dbad}, we similarly choose interpolation elements within $\Omega=[0,2\pi)^2$ according to different target points. 
\Cref{tab:example4_fpr_error} shows that the recovery error gradually decreases as the interpolation elements decrease in size. 
Note that the interpolated element size cannot be exactly halved as in the three examples above due to the selection of the computable region $\mathcal{G}$, leading to the accuracy order of the FPR method exhibits fluctuations. 

\begin{table}[!bpth]
    \centering
        \caption{Error of FPR method with 2D LIF-1 of recovering quasiperiodic function \cref{eq:approx} with interpolation element sizes $\langle\Theta\rangle$.}
    \label{tab:example4_fpr_error}
    \begin{tabular}{|c|c|c|}
    \hline
    $\langle\Theta\rangle$  & Error & Order\\
	\hline
	(0.4048, 0.3)  & 1.1717e-01 & \\
	(0.2024, 0.15)  & 2.0538e-02 & 2.51\\
	(0.1012, 0.075)  & 5.1695e-03 & 1.99\\
	(0.0667, 0.0375) & 1.6952e-03 & 1.61\\
	(0.0328, 0.0188) & 4.2629e-04 & 1.99\\
    \hline
    \end{tabular}
\end{table}

\subsection{Piecewise constant quasicrystals}

\begin{figure*}[!hbtp]
    \centering
	\subfigure[1D Fibonacci photonic quasicrystal structure in 2D tilling]{\label{fig:fibonacci1}\includegraphics[width=6.5cm]{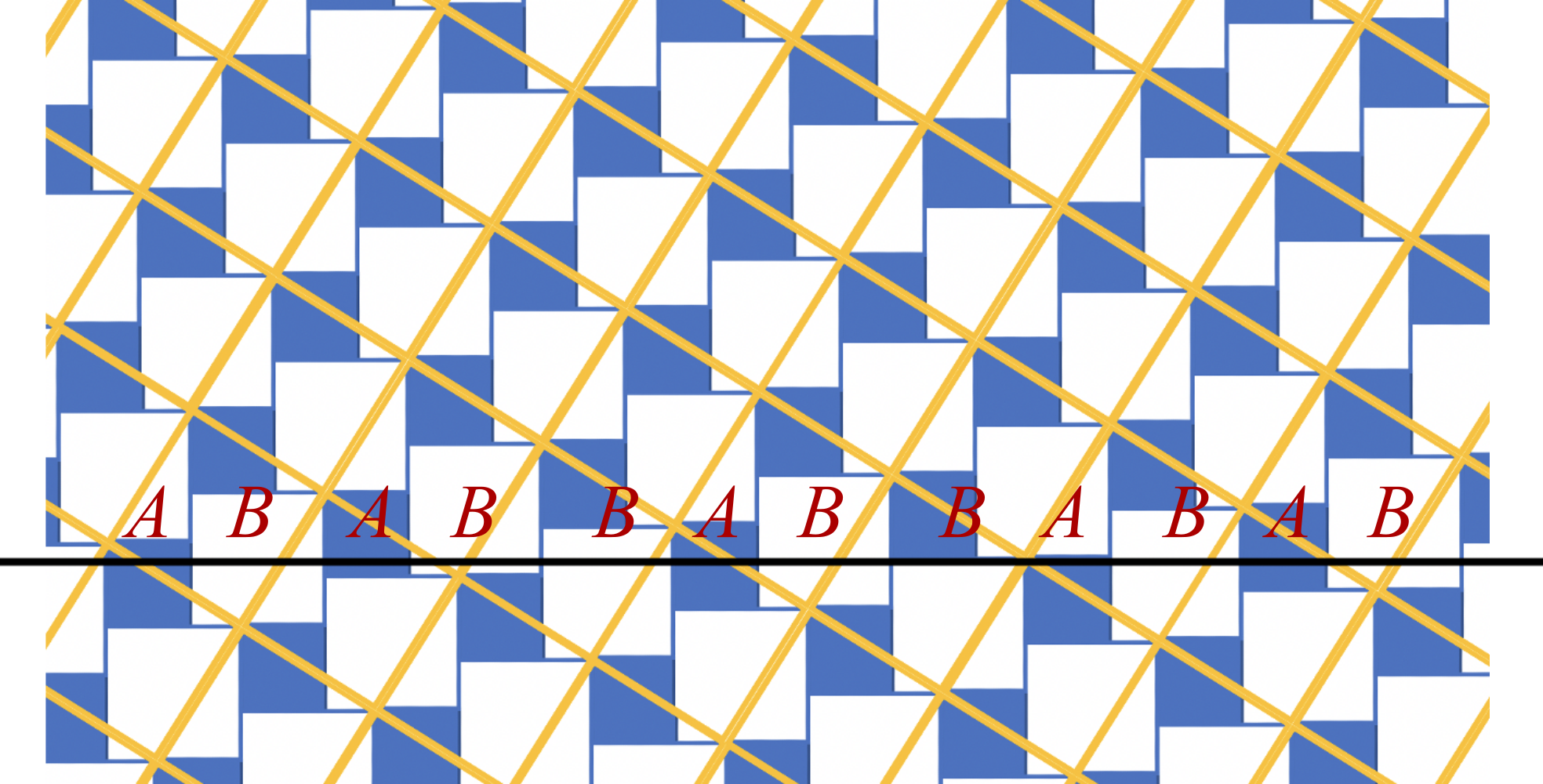}}
	\subfigure[Parent function $F(\bm{s})$ in $\Omega$]{\label{fig:fibonacci3}\includegraphics[width=4.5cm]{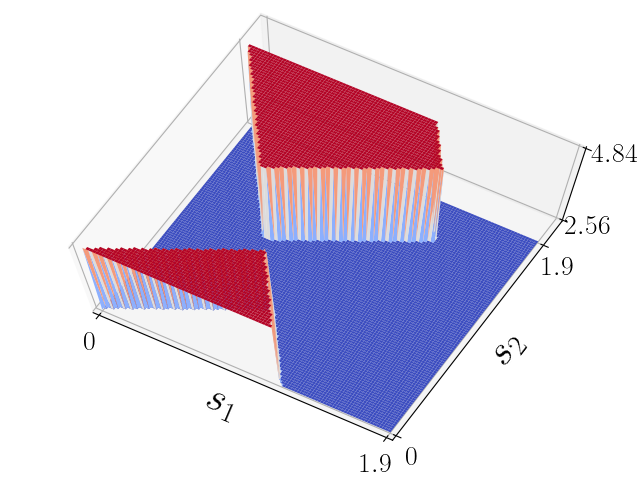}}
	\caption{1D Fibonacci photonic quasicrystal and its parent function $F(\bm{s})$ in $\Omega=[0,1.9)^2$.}
\end{figure*}

\begin{example}\label{exa:fibonacci}
    Consider a 1D Fibonacci photonic quasicrystal as shown in \Cref{fig:fibonacci1}, which is a black line in a 2D tilling~\cite{vardeny2013optics, tanese2014fractal}. 
    The blue (white) square has side length $A$ ($B$), and its corresponding dielectric constant is $\varepsilon_A=4.84$ ($\varepsilon_B=2.56$). 
    The sequence satisfies $F_n=F_{n-2}+F_{n-1},~n=1,2,\cdots$, with $F_1=B$,  $F_2=A$, and  
    $\lim\limits_{j\to\infty}F_{j+1}/F_{j}=(1+\sqrt{5})/2:=\lambda$.
\end{example}

Let the angle $\phi$ be defined by $\tan \phi =\lambda$. 
The 2D tilling can be transformed into a periodic structure by a rotation matrix defined as
\begin{equation*}
    \begin{bmatrix}
        \sin \phi & \cos \phi \\
        \cos \phi & -\sin \phi \\
    \end{bmatrix}.
\end{equation*}
All unit cells are surrounded by the yellow border in \Cref{fig:fibonacci1}. 
\Cref{fig:fibonacci3} shows the dielectric constant over the unit cell $\Omega=[0,1.9)^2$. 
Therefore, the projection matrix of 
the 1D Fibonacci photonic quasicrystal is $\bm{P}=(\sin \phi, \cos \phi)$.

\begin{figure*}[!hbpt]
    \centering
	\subfigure[Exact value]{\label{fig:fibonacci5_1}\includegraphics[width=5cm]{./figures/fibonacci_exact.png}}
	\hspace{10mm}
	\subfigure[FPR with 2D LIF-1]{\label{fig:fibonacci5_2}\includegraphics[width=5cm]{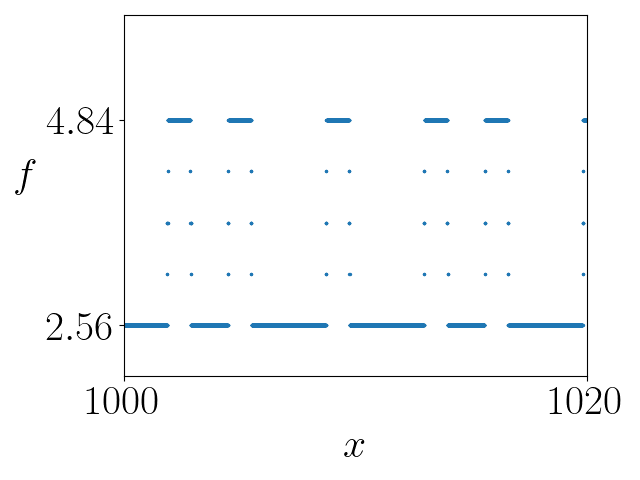}}
	\caption{Comparison of the exact value and the numerical result obtained by FPR method for the dielectric constant of 1D Fibonacci photonic quasicrystal when $x\in[1000,1020)$.}
    \label{fig:fibonacci5}
\end{figure*}

\begin{figure*}[!hbpt]
    \centering
	\includegraphics[width=0.8\textwidth]{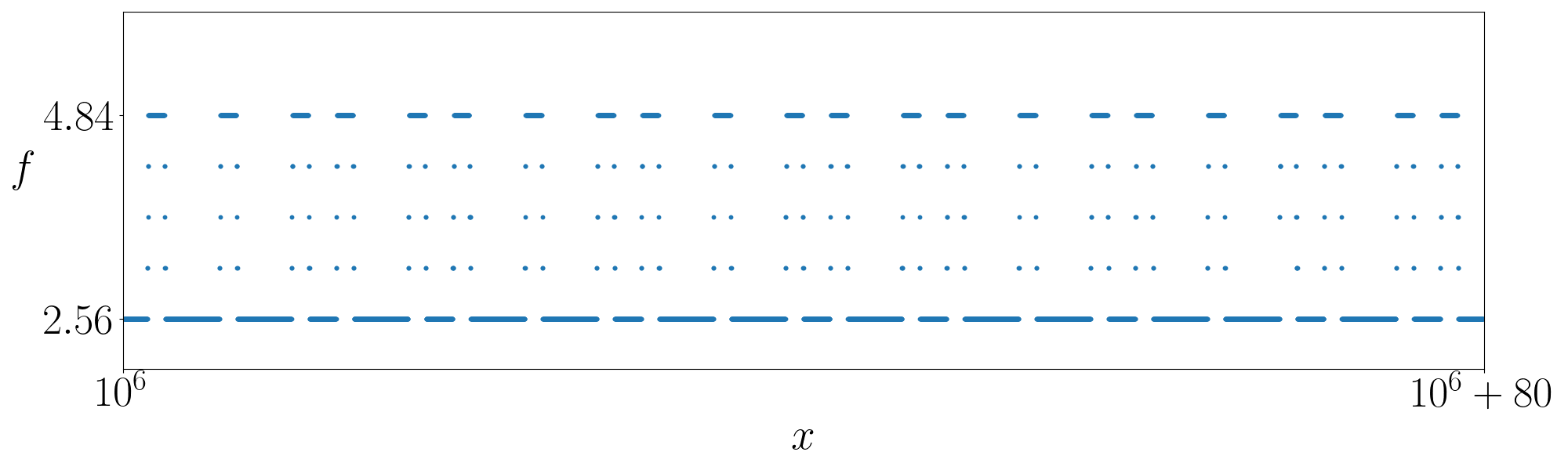}
	\caption{Dielectric constant of 1D Fibonacci photonic quasicrystal recovered by FPR method with 2D LIF-1 when $x\in[10^6,10^6+80)$.}
    \label{fig:fibonacci_global}
\end{figure*}

Given the filling precision $\varepsilon=10^{-3}$, we select the computable region $\mathcal{G}=[0,500)$ according to \Cref{tab:suboptimal_region}. 
When we select the size of interpolation elements $\langle\Theta\rangle=(0.08,0.08)$, only 480 interpolation nodes are needed to recover the global dielectric constant of the 1D Fibonacci photonic quasicrystal. 
\Cref{fig:fibonacci5_2} shows the recovered dielectric constant of the 1D Fibonacci photonic quasicrystal when $x\in[1000,1020)$. 
FPR method can accurately recover the dielectric constant at continuous point compared with the exact value as shown in \Cref{fig:fibonacci5_1}.
Moreover, FPR method is capable to seize the position of discontinuity points. 
Meanwhile, we present \Cref{fig:fibonacci_global} to demonstrate the effectiveness of the FPR method for recovering global function.

\section{Conclusion and outlook}\label{sec:con}
This paper is concerned with developing a new algorithm for recovering both smooth and non-smooth quasiperiodic systems. 
Different from the existing spectral Galerkin methods, we propose the FPR method for accurately recovering the global quasiperiodic system using interpolation technique based on finite points. 
To theoretically support our method, we establish a homomorphism between the physical space of quasiperiodic function and the high-dimensional torus. 
Moreover, by exploiting the arithmetic properties of irrational numbers, we design the FPR method with exquisite algorithmic steps to ensure accuracy and efficiency in the recovery process. 
We also present the corresponding convergence analysis and computational complexity analysis. 
We apply our algorithm to solve two classes of quasiperiodic problems: continuous quasiperiodic functions and a piecewise constant Fibonacci  quasicrystal. 
Numerical results show the effectiveness and superiority of FPR approach, while PM method fails to recover the non-smooth quasiperiodic systems.
Furthermore, the experiments conclusively demonstrate that the FPR method, as a non-lifting algorithm, exhibits substantial computational advantages compared to the state-of-the-art projection method. 

There remains much work to be done based on the proposed method. 
First, we aim to extend the convergence analysis of FPR method to a more general form that can accommodate non-smooth parent functions. 
Second, we intend to introduce a parallel strategy and improve the selection method of sampling points in the algorithm implementation, especially for high-dimensional problems. 
Third, we plan to develop the FPR method to handle singularly quasiperiodic systems by integrating the adaptive strategy. 
Finally, we will apply the FPR method to solve physical problems, discover exotic phenomena and physical law. 
These above-mentioned endeavors will not only enrich the theoretical foundations of FPR method, but also contribute to practical applications in diverse fields, such as physics, engineering, and materials science.

% \bibliographystyle{siamplain}
% \bibliography{references.bib}

\end{document}